\documentstyle[12pt]{amsart}
\newcommand{\g}{\goth}
\newcommand{\gtg}{\mbox{\g g}}

\newcommand{\hgtg}{\mbox{$\hat{\gtg}$}}

\newcommand{\gtsl}{\mbox{\g sl}}
\newcommand{\gtgl}{\mbox{\g gl}}
\newcommand{\hgtsl}{\mbox{$\hat{\gtsl}$}}
\newcommand{\hgtgl}{\mbox{$\hat{\gtgl}$}}

\newcommand{\gth}{\mbox{\g h}}
\newcommand{\gtS}{\mbox{\g S}}
\newcommand{\nc}{\mbox{${\Bbb C}$}}
\newcommand{\nz}{\mbox{${\Bbb Z}$}}

\newcommand{\cA}{\mbox{${\cal A}$}}

\newcommand{\cD}{\mbox{${\cal D}$}}
\newcommand{\cF}{\mbox{${\cal F}$}}

\newcommand{\cP}{\mbox{${\cal P}$}}
\newcommand{\cR}{\mbox{${\cal R}$}}
\newcommand{\cN}{\mbox{${\cal N}$}}
\newcommand{\id}{\mbox{\rm id}}
\newcommand{\qdet}{\mbox{\rm q-det}}
\newcommand{\sgn}{\mbox{\rm sgn}}
\newcommand{\End}{\mbox{\rm End}}

\newcommand{\rank}{\mbox{\rm rank}}
\newcommand{\Ad}{\mbox{\rm Ad}}
\newcommand{\vep}{\varepsilon}

\theoremstyle{plain}
 \newtheorem{thm}{Theorem}[section]
 \newtheorem{prop}[thm]{Proposition}
 \newtheorem{lemma}[thm]{Lemma}
 \newtheorem{cor}[thm]{Corollary}
\theoremstyle{definition}
 \newtheorem{defn}{Definition}[section]
 \newtheorem{conj}{Conjecture}
 \newtheorem{ex}{Example}[section]
\theoremstyle{remark}
 \newtheorem{rem}{Remark}

 \newtheorem{ack}{Acknowledgement}

\begin{document}
\title{Bosonic representations of Yangian Double $\cD Y_{\hbar}(\gtg)$  
               with $\gtg=\gtgl_N,\gtsl_N$                   \\}
\author{Kenji Iohara$^\dagger$}
\address[K. Iohara]{Department of Mathematics, Faculty of Science,
                    Kyoto University, Kyoto 606, Japan.}
\thanks{$\dagger$:JSPS Research Fellow}
\date{\today}
\maketitle
\begin{abstract}
On the basis of `$RTT=TTR$' formalism, we introduce the quantum double 
of the Yangian $Y_{\hbar}(\gtg)$ for $\gtg=\gtgl_N,\gtsl_N$ with a central
extension. The Gauss decomposition of T-matrices gives us the so-called
Drinfel'd generators. Using these generators, we present some examples
of both finite and infinite dimensional representations that are quite
natural deformations of the corresponding affine counterpart.
\end{abstract}
\pagestyle{plain}
\section{Introduction}\hspace{1 in} \\

For the last few decades, the quantum inverse scattering method (QISM), 
initiated by L. D. Faddeev and his colleagues, have been studied extensively 
and produced rich structures both in physics and in mathematics. Quantum
algebras called quantized universal enveloping algebra $U_q(\gtg)$ and
Yangian $Y_{\hbar}(\gtg)$ are one of the most important fruits inspired
by the QISM. They have unexpected connections with such, at first sight, 
unrelated parts of mathematics as the construction of knot invariants, the
geometric interpretation of a certain class of special functions and the
representation theory of algebraic groups in characteristic $p$. Of course
they also have many nice applications in such theoretical physics as 
quantum field theory and statistical mechanics. As is well-known, 
$U_q(\gtg)$ describes some features of conformal field theory. One can
solve lattice models, like the spin $\frac{1}{2}$ XXZ model, as an
application of the representation theory of $U_q(\hgtsl_2)$ . 
The quantum affine algebra $U_q(\hgtsl_2)$ is the 
$q$-deformation of the enveloping algebra $U(\hgtsl_2)$. 
The Yangian $Y_{\hbar}(\gtg)$ is also related to conformal field theory.
Lattice models such as the Haldane-Shastry model are known to possess 
$Y_{\hbar}(\gtsl_2)$-symmetry. The Yangian $Y_{\hbar}(\gtsl_2)$ is the 
$\hbar$-deformation of the enveloping algebra $U(\gtsl_2[t])$.
The quantum double \cite{dr_icm} of the $Y_{\hbar}(\gtg)$, 
which we shall refer to as Yangian double $\cD Y_{\hbar}(\gtg)$, 
seems to play important roles in massive field theory 
\cite{bern},\cite{leclair},\cite{smirnov}. In these works, the Yangian double
$\cD Y_{\hbar}(\gtg)$ is the $\hbar$-deformation of the universal 
enveloping algebra of the loop algebra $\gtg[t,t^{-1}]$ for $\gtg =\gtsl_2$, 
{\it without central extension}. In view of lattice models, like the
spin $\frac{1}{2}$ XXX model of infinite chains, it seems necessary to 
construct the Yangian Double $\cD Y_{\hbar}(\gtg)$ with a central extension.
In our previous paper \cite{kohno}, we defined the Yangian double 
$\cD Y_{\hbar}(\gtg)$ with a central extension for $\gtg=\gtgl_2$ or 
$\gtsl_2$. The present paper is a higher rank generalization of it. 
Our attempt here is to explain the background of the construction and 
to consider the representation theory. We also summarize some formulae 
related to our calculations which seem well-known to the specialists 
but have never appeared in the literature. The main topics treated 
in this paper is as follows.

\noindent{\bf 1. Yangian Double}
Yangian double $\cD Y_{\hbar}(\gtg)$ has been introduced in the literature 
in terms of Chevalley generators \cite{leclair}, $T^{\pm}$-matrix \cite{bern}
for $\gtg =\gtsl_2$ and Drinfel'd generators \cite{tolstoy} for a simple finite
dimensional Lie algebra $\gtg$. Here we construct $\cD Y_{\hbar}(\gtg)$ for
$\gtg =\gtgl_N,\gtsl_N$ by means of the QISM \cite{bern}, \cite{resh},
\cite{QISM}. Namely, let $R(u)$ be the Yang's $R$-matrix. The algebra 
$\cD Y_{\hbar}(\gtgl_N)$ is defined through quadratic relations of the form
\begin{align*}
R(u-v)(T^{\pm}(u)\otimes \id)(\id \otimes T^{\pm}(v))&=
(\id \otimes T^{\pm}(v))(T^{\pm}(u)\otimes \id)R(u-v), \\
R(u-v-\frac{1}{2}\hbar c)(T^+(u)\otimes \id)(\id \otimes T^-(v))&=
(\id \otimes T^-(v))(T^+(u)\otimes \id)R(u-v+\frac{1}{2}\hbar c),
\end{align*} 
where $c$ is a central element of $\cD Y_{\hbar}(\gtgl_N)$. The 
$T^{\pm}$-matrix $T^{\pm}(u)=(t_{ij}^{\pm}(u))_{1\leq i,j \leq N}$
are expanded as
\[ t_{ij}^+(u)=\delta_{ij}-\hbar \sum_{k\geq 0} t_{ij}^k u^{-k-1}, \quad
   t_{ij}^-(u)=\delta_{ij}+\hbar \sum_{k<0} t_{ij}^k u^{-k-1}. \]
Just as in the case of $U_q(\hgtgl_n)$ \cite{ding}, we consider the Gauss
decomposition of $T^{\pm}$-matrix (Theorem \ref{ding_frenkel}) and obtain
the Drinfel'd generators of $\cD Y_{\hbar}(\gtgl_N)$ (Theorem \ref{gl_N}).
We define $\cD Y_{\hbar}(\gtsl_N)$ as a certain subalgebra of 
$\cD Y_{\hbar}(\gtgl_N)$ and show that our Drinfel'd generators recover
the results obtained in \cite{tolstoy} at level $0$ (Corollary \ref{comm_0}).
We also introduce another subalgebra of $\cD Y_{\hbar}(\gtgl_N)$ 
which we call Heisenberg subalgebra.

\noindent{\bf 2. Representation Theory}
Here we investigate several examples. The main tool here is Drinfel'd
generators.

\noindent{\bf Finite dimensional representations}
At $c=0$, the Heisenberg subalgebra of $\cD Y_{\hbar}(\gtgl_N)$ becomes
the center of it. So we will concentrate on $\cD Y_{\hbar}(\gtsl_N)$ case
without loss of generality. From the commutation relations of 
$\cD Y_{\hbar}(\gtg)$ at level $0$ (Corollary \ref{comm_0}), 
we expect that the analogue of the classification 
theorem of irreducible finite dimensional representations holds 
just as in the case of Yangian $Y_{\hbar}(\gtg)$ \cite{dr_2}. We present some
examples which support our conjecture. All of them are those what we call
evaluation modules.

\noindent{\bf Infinite dimensional representations}
Unfortunately we have no proper definition of highest weight modules due
to the lack of the triangular decomposition of $\cD Y_{\hbar}(\gtg)$. 
Here we realize level $1$ $\cD Y_{\hbar}(\gtgl_N)$-modules on boson Fock
space $\cF_{i,s}$~($0\leq i\leq N-1, s\in \nc$) (Theorem \ref{level 1 gl_N}). 
Let $V_u$ be an $N$-dimensional evaluation modules of $\cD Y_{\hbar}(\gtgl_N)$.
Vertex operators are the intertwiners of the form
\begin{align*}
\Phi^{(i,i+1)}(u) & : \cF_{i+1,s}\longrightarrow
                      \cF_{i,s-1}\otimes V_{u},  \\
\Psi^{(i,i+1)}(u) & : \cF_{i+1,s}\longrightarrow
                      V_{u}\otimes \cF_{i,s-1}. 
\end{align*}
We also give the bosonization of vertex operators (Theorem \ref{VOp_gl_N}).
For the $\cD Y_{\hbar}(\gtsl_N)$ case, we construct level $1$ modules on
boson Fock space $\cF_i$~($0\leq i\leq N-1$) (Theorem \ref{level 1 sl_N})
whose quantum affine version are obtained in \cite{jing}. We should
mention that every field defined above makes sense as a formal series in
$\hbar$. Moreover, we also construct vertex operators for 
$\cD Y_{\hbar}(\gtsl_N)$, in which case the Fourier components lose
its meaning (Theorem \ref{VOp_sl_N}). More precisely, those formulae makes
sense only as an asymptotic series.

 The text is organized as follows. In Section 2 we recall the definition
of Yangian $Y_{\hbar}(\gtg)$. We also mention about the other set of 
generators and the isomorphism among them. Theory of finite dimensional 
$Y_{\hbar}(\gtg)$-modules is also reviewed and one example is given.
In Section 3 we define $\cD Y_{\hbar}(\gtg)$ for $\gtg=\gtgl_N,\gtsl_N$.
We rewrite the commutation relations in terms of Drinfel'd generators.
In Section 4 we present a conjecture for finite dimensional 
$\cD Y_{\hbar}(\gtsl_N)$-modules together with a few examples. 
As for infinite dimensional representations, we construct level 1 modules
and vertex operators directly via bosonization. Section 5 contains 
discussions and remarks. For the reader's convenience, we also include
two appendices. In Appendix A we give a brief review of the quantum group 
especially about universal $\cR$ and $L$-operators. 
In Appendix B we collect some formulas of $T$-matrices.

 Let us mention that the author got two papers \cite{khoroshkin1},
\cite{khoroshkin2} when he was preparing this article. The central extension
of $\cD Y_{\hbar}(\gtsl_2)$ is introduced in \cite{khoroshkin1} 
which has some overlap with \cite{kohno}. The bosoonizations of level $1$ 
$\cD Y_{\hbar}(\gtsl_2)$-module and the vertex operators among them are
obtained in \cite{khoroshkin2}. Here we introduce the Yangian Double
$\cD Y_{\hbar}(\gtg)$ for $\gtg=\gtgl_N,\gtsl_N$ with a center and obtain
the bosonization of level $1$ $\cD Y_{\hbar}(\gtg)$-module and the vertex 
operators among them.  
\section{Review of Yangian $Y_{\hbar}(\gtg)$}\hspace{1 in} \\
In this section we collect some known facts about Yangians including 
representation theory.

\subsection{Yangian $Y_{\hbar}(\gtg)$}
Here we present two different realizations of $Y_{\hbar}(\gtg)$ for
a simple finite dimensional Lie algebra $\gtg$. In addition, for
$\gtg=\gtsl_N$, another realization called $T$-matrix is known
 \cite{dr_icm,dr_1}, and we give some comments on it.

Set $\cA=\nc[[\hbar]]$. Let $\gtg$ be a simple finite 
dimensional Lie algebra and $\{\alpha_1,\alpha_2,\cdots ,\alpha_n \}$ the set
of simple roots. Fix a standard non-degenerate symmetric invariant
bilinear form $(\cdot,\cdot)$ on $\gtg$. For each positive root $\alpha$ of
$\gtg$, choose root vectors $x^{\pm}_{\alpha}$ in $\pm \alpha$ root spaces
such that $(x^{+}_{\alpha},x^{-}_{\alpha})=1$ and set 
$h_{\alpha}=[x^{+}_{\alpha},x^{-}_{\alpha}]$. We denote the Cartan matrix of
$\gtg$ by $A=(a_{ij})$.
 
Let $\{I_p \}$ be any orthonormal basis of $\gtg$ with respect to 
the inner product $(\cdot, \cdot)$.
\begin{defn}[\cite{dr_2}]
The Yangian $Y_{\hbar}(\gtg)$ is a topological Hopf algebra over $\cA$
generated by $\gtg$ and elements $J(x),x\in \gtg$, with relations
\begin{align*}
& J(ax+by)=aJ(x)+bJ(y),~a,b \in \cA, [x,J(y)]=J([x,y]), \\
& [J(x),J([y,z])]+[J(y),J([z,x])]+[J(z),J([x,y])]\\
& =\hbar^2 \sum_{p,q,r}([x,I_p],[[y,I_q],[z,I_r]])\{I_p,I_q,I_r\}, \\
& [[J(x),J(y)],[z,J(w)]]+[[J(z),J(w)],[x,J(y)]]\\
& =\hbar^2 \sum_{p,q,r}
   \left( ([x,I_p],[[y,I_q],[[z,w],I_r]])+
          ([z,I_p],[[w,I_q],[[x,y],I_r]]) \right) \{I_p,I_r,J(I_r) \}, 
\end{align*}
where $\{\cdot, \cdot, \cdot \}$ denotes the symmetrization
\[ \{ x_1,x_2,x_3 \}=\frac{1}{24}\sum_{\sigma \in \gtS_3}
   x_{\sigma(1)}x_{\sigma(2)}x_{\sigma(3)}. \]
The comultiplication of $Y_\hbar(\gtg)$ is given by
\begin{align*}
& \Delta(x)=x\otimes 1 + 1 \otimes x, \\
& \Delta(J(x))=J(x)\otimes 1 + 1 \otimes J(x) 
  + \frac{1}{2}\hbar[x\otimes 1, \Omega], 
\end{align*}
where $\Omega$ stands for the Casimir element of $\gtg \otimes \gtg$.
\end{defn}
It is known by Drinfel'd \cite{dr_2} that there is so-called Drinfel'd
generators of $Y_\hbar(\gtg)$. To be precise, the following theorem holds.
\begin{thm}[\cite{dr_2}]
The Yangian $Y_{\hbar}(\gtg)$ is isomorphic to the algebra
generated by the elements 
$\{\xi_{ik}^{\pm},\kappa_{ik}
  |1\leq i\leq n~,~k\in \nz_{\geq 0} \}$ subject to the relations
\begin{align*}
& [\kappa_{ik},\kappa_{jl}]=0~,
 ~[\kappa_{i0},\xi_{jl}^{\pm}]=\pm(\alpha_i,\alpha_j)\xi_{jl}^{\pm}~,
 ~[\xi_{ik}^+,\xi_{jl}^-]=\delta_{ij}\kappa_{ik+l}~, \\ 
& [\kappa_{ik+1},\xi_{jl}^{\pm}]-[\kappa_{ik},\xi_{jl+1}^{\pm}]
 =\pm \frac{1}{2}(\alpha_i,\alpha_j)\hbar 
  [\kappa_{ik},\xi_{jl}^{\pm}]_+~, \\
& [\xi_{ik+1}^{\pm},\xi_{jl}^{\pm}]
 -[\xi_{ik}^{\pm},\xi_{jl+1}^{\pm}]
 =\pm \frac{1}{2}(\alpha_i,\alpha_j)\hbar
  [\xi_{ik}^{\pm},\xi_{jl}^{\pm}]_+~, \\
& \sum_{\sigma \in \gtS_{m}}
  [\xi_{ik_{\sigma(1)}}^{\pm},[\cdots,
  [\xi_{ik_{\sigma(m)}}^{\pm},\xi_{jl}^{\pm}]=0,\quad \text{for}
  \quad i\neq j~,
\end{align*}
where we set $m=1-a_{ij}$ and $[x,y]_+=xy+yx$ for 
             $x,y \in Y_{\hbar}(\gtg)$. 
The isomorphism $\phi$ between two presentations is given by
\begin{align*}
 \phi(h_i)=\kappa_{i0} &, \quad \phi(x^{\pm}_i)=\xi^{\pm}_{i0}, \\
 \phi(J(h_i))=\kappa_{i1}+\hbar \phi(v_i) &, \quad
 \phi(J(x^{\pm}_i))=\xi^{\pm}_{i1}+\hbar \phi(w^{\pm}_i),
\end{align*}
where we set $h_i=h_{\alpha_i},x^{\pm}_i=x^{\pm}_{\alpha_i}$ and 
\begin{align*}
& v_i=\frac{1}{4}\sum_{\alpha \succ 0}(\alpha,\alpha_i)
                 (x^{+}_{\alpha}x^{-}_{\alpha}+x^{-}_{\alpha}x^{+}_{\alpha})
     -\frac{1}{2}{h_i}^2, \\
& w^{\pm}_i=\pm \frac{1}{4}\sum_{\alpha \succ 0}
            \{ [x^{\pm}_i,x^{\pm}_{\alpha}]x^{\mp}_{\alpha}+
               x^{\mp}_{\alpha} [x^{\pm}_i,x^{\pm}_{\alpha}] \}
           -\frac{1}{4}(x^{\pm}_i h_i + h_i x^{\pm}_i).
\end{align*}
\end{thm}
For $\gtg=\gtsl_N$, we have another 
realization called $T$-matrix \cite{dr_1,dr_2} as follows.

Let $V$ be a rank $N$ $\cA$-free module
and  $\cP \in \End(V\otimes V)$
be a permutation operator $\cP v\otimes w =w\otimes v~(v,w\in V)$.
Consider Yang's $R$-matrix normalized as
\begin{equation}\label{Yang's R}
 R(u)=\frac{1}{u+\hbar}\left( uI+\hbar\cP \right) \in
        \End(V\otimes V).
\end{equation}
where $\hbar$ is expanded in positive powers.
This $R$-matrix satisfies the following properties:
\begin{description}
\item[Yang-Baxter equation]
\[ R_{12}(u-v)R_{13}(u)R_{23}(v)=
   R_{23}(v)R_{13}(u)R_{12}(u-v), \]
\item[Unitarity]
\[ R_{12}(u)R_{21}(-u)=\id . \]
\end{description}
Here, if $R(u)=\sum a_i \otimes b_i$ with $a_i,b_i \in \End (V)$, then
$R_{21}(u)=\sum b_i \otimes a_i,~R_{13}(u)=\sum a_i \otimes 1 \otimes b_i$
etc..

\begin{thm}[\cite{dr_2}]
The Yangian $Y_{\hbar}(\gtsl_N)$ is isomorphic to the algebra with 
generators $\{ t_{i,j}^k | 1\leq i,j\leq N~,~k \in \nz_{\geq 0} \}$
and defining relations
\[ R(u-v)\overset{1}{T}(u)\overset{2}{T}(v)=
   \overset{2}{T}(v)\overset{1}{T}(u)R(u-v),\quad \qdet. T(u)=1. \]
Here 
\begin{align*}
& T(u)=(t_{ij}(u))_{1\leq i,j \leq N},\quad 
   t_{ij}(u)=\delta_{ij}
             -\hbar \sum_{k\in \nz_{\geq 0}}t_{ij}^{k}u^{-k-1}, \\
& \overset{1}{T}(u)=T(u)\otimes \id~,~\overset{2}{T}(u)=\id \otimes T(u),
\end{align*}
and $\qdet. T(u)$ is defined in $($\ref{defn_q-det}$)$. 
The comultiplication is given by
\[ \Delta(t_{ij}(u))=\sum_{k=1}^N t_{kj}(u)\otimes t_{ik}(u). \]
\end{thm}
Roughly speaking, the isomorphism between the algebra generated by the 
Drinfel'd generators and the algebra presented above is given by the 
Gauss decomposition of the $T$-matrix (See Section \ref{main} 
and Appendix \ref{gauss} for detail.).
\subsection{Representation theory of $Y_{\hbar}(\gtg)$}
In this subsection, we give a brief review on finite dimensional
representations of $Y_{\hbar}(\gtg)$. See \cite{chari1,chari2} for detail.

Let {\bf h}$=\{ h_{i,r} \}_{1\leq i \leq n, r\in \nz_{\geq 0}}$ be a 
subset of $\cA$. A $Y_{\hbar}(\gtg)$-module $V$ is called {\it highest weight
module with highest weight {\bf h}} if there exits an unique, up to scalar,
non-zero vector $v\in V$ such that $V$ is generated by $v$ and 
\[ \kappa_{i,r}.v=h_{i,r}v,\quad \xi^{+}_{i,r}.v=0, \qquad 
   1\leq \forall i \leq n,~\forall r \in \nz_{\geq 0}. \]
It is known that every irreducible finite dimensional 
$Y_{\hbar}(\gtg)$-module $V$ is highest weight module. Let us denote the 
irreducible highest weight $Y_{\hbar}(\gtg)$-module with highest weight
{\bf h} by $V(\text{{\bf h}})$. The criterion of the finite dimensionality
of $V(\text{{\bf h}})$ is known.
\begin{thm}[\cite{dr_2}]
The irreducible $Y_{\hbar}(\gtg)$-module $V(\text{{\bf h}})$ of highest
weight {\bf h} is finite dimensional if and only if there exist monic
polynomials $P_{i}(v) \in \cA [v]~1\leq i \leq n$ such that
\[ \frac{P_{i}(v+\frac{1}{2}(\alpha_{i},\alpha_{i})\hbar)}{P_{i}(v)}=
   1+\hbar \sum_{r=0}^{\infty}h_{i,r}v^{-r-1}, \]
in the sense that the right-hand side is the Laurent expansion of the
left-hand side about $v=\infty$.
\end{thm}
The Polynomials $P_{i}(v)$ in this theorem are called Drinfel'd polynomials.

The preceding theorem suggests the following definition.
\begin{defn}[\cite{chari2}]
We say that an irreducible finite dimensional $Y_{\hbar}(\gtg)$-module is
{\it fundamental} if its Drinfel'd polynomials are given by
\[ P_j(v)=\begin{cases} v-u \quad j=i \\
                         \quad 1  \qquad j\neq i
          \end{cases} \]
for some $1 \leq i \leq n$.
\end{defn}
Using the fact that the Drinfel'd polynomials of the tensor product of two 
highest weight modules are the product of the Drinfel'd polynomials of
two highest weight modules, one proves the following.
\begin{thm}[\cite{chari2}]
Every irreducible finite dimensional $Y_{\hbar}(\gtg)$-module is isomorphic
to a subquotient of a tensor product of fundamental representations.
\end{thm}
 
Here we present an example of fundamental representation of 
$Y_{\hbar}(\gtsl_N)$ and the more general representations containing 
the first example. 


\begin{ex}\label{Y(sl_N)}
 Set
\[ V_{u}=V\otimes_{\cA} \cA[u],\quad V=\oplus_{j=0}^{N-1}\cA w_j. \]
$Y_{\hbar}(\gtsl_N)$-module structure on $V_{u}$ is defined via the following
actions.
\begin{align*}
(1)\qquad & \xi^{+}_{i,r}.w_{i}=(u-\frac{N-1-i}{2}\hbar)^r.w_{i-1}, \qquad
            \xi^{+}_{i,r}.w_{j}=0 \quad j\neq i, \\
(2)\qquad & \xi^{-}_{i,r}.w_{i-1}=(u-\frac{N-1-i}{2}\hbar)^r.w_{i}, \qquad
            \xi^{-}_{i,r}.w_{j}=0 \quad j\neq i-1,  \\
(3)\qquad & \kappa_{i,r}.w_{i-1}=(u-\frac{N-1-i}{2}\hbar)^r.w_{i-1}, \qquad
            \kappa_{i,r}.w_{i}=-(u-\frac{N-1-i}{2}\hbar)^r.w_{i}, \\
          & \kappa_{i,r}.w_{j}=0 \quad j\neq i,i-1. 
\end{align*}
Note that one can regard $u$ as either an indeterminate or an element of $\cA$.
In the latter case, the Drinfel'd polynomials of $V_{u}$ are given by
\[ P_{1}(v)=v-(u-\frac{N-2}{2}\hbar),\qquad P_{i}(v)=1 \quad i\neq 1. \]
\end{ex}

Let us fix $\gtg$ to be a simple finite dimensional Lie algebra of 
classical type and normalize the invariant bilinear form by the condition
$(\beta,\beta)=2$~($\beta$: long root). The fundamental weights 
$\Lambda_i (1\leq i\leq \rank \gtg)$ are chosen so as to satisfy
\[ 2\frac{(\Lambda_i,\alpha_j)}{(\alpha_j,\alpha_j)}=\delta_{i,j}. \]

\begin{ex}\label{generic_repr} 
Let $V(\Lambda)$ be the irreducible highest weight $\gtg$-module
with highest weight $\Lambda$. Especially when $\Lambda$ is of the form
$m\Lambda_i$ with $m$ being positive integer, it is known by \cite{kir}
if $\gtg$ is of type $A_l$(resp. $B_l,C_l,D_l$) and $1\leq i \leq l$
(resp. $i=1,i=l,i=1,l-1,l$) then $V(m\Lambda_i)$ can be made into 
$Y_{\hbar}(\gtg)$-module. Let $V_a(m\Lambda_i)$ be such 
$Y_{\hbar}(\gtg)$-module satisfying
\begin{align*}
(1)\qquad & V_a(m\Lambda_i)\cong V(m\Lambda_i)\quad \text{as $\gtg$-module}, \\
(2)\qquad & J(x)\vert_{V_a(m\Lambda_i)}=ax\vert_{V_a(m\Lambda_i)}\quad
            \forall x \in \gtg. 
\end{align*}
The Drinfel'd polynomials of $V_a(m\Lambda_i)$ are given by
\[ P_{i}(v)=\prod_{k=1}^m\{ u-(a-(\frac{1}{4}g-j+\frac{m}{2})\hbar)\},
   \qquad P_{j}(v)=1 \quad j\neq i, \]
where $g$ is the dual Coxeter number.
\end{ex}
\section{The algebra $\cD Y_{\hbar}(\gtgl_N)$}\label{main}\hspace{1 in} \\ 
Here we define a central extension of $\cD Y_{\hbar}(\gtg)$ for 
$\gtg=\gtgl_N,\gtsl_N$ following the method of \cite{resh}.

\subsection{Yangian Double $\cD Y_{\hbar}(\gtgl_N)$}
Let us choose Yang's $R$-matrix as in (\ref{Yang's R}).

\begin{defn}\label{defn_DY(gl_N)}
$\cD Y_{\hbar}(\gtgl_N)$ is a topological Hopf algebra over
$\cA$ generated by
$\{ t_{ij}^k |1\leq i,j \leq N~,~k\in \nz \}$ and $c$. In terms of matrix
generating series
\begin{align*}
&  T^{\pm}(u)=(t_{ij}^{\pm}(u))_{1\leq i,j \leq N}, \\
&  t_{ij}^{+}(u)=\delta_{ij}
                -\hbar \sum_{k\in \nz_{\geq 0}}t_{ij}^{k}u^{-k-1}~,~
   t_{ij}^{-}(u)=\delta_{ij}
                +\hbar \sum_{k\in \nz_{<0}}t_{ij}^{k}u^{-k-1}, 
\end{align*}
the defining relations are given as follows:
\begin{align*}
[T^{\pm}(u),c]& =0, \\
R(u-v)\overset{1}{T^{\pm}}(u)\overset{2}{T^{\pm}}(v)&=
\overset{2}{T^{\pm}}(v)\overset{1}{T^{\pm}}(u)R(u-v), \\
R(u_{-}-v_{+})\overset{1}{T^{+}}(u)\overset{2}{T^{-}}(v)&=
\overset{2}{T^{-}}(v)\overset{1}{T^{+}}(u)R(u_{+}-v_{-}).
\end{align*}
Here
\[   \overset{1}{T}(u)=T(u)\otimes \id~,
    ~\overset{2}{T}(u)=\id \otimes T(u), \]
$u_{\pm}=u\pm \frac{1}{4}\hbar c$ and similarly for $v$.
Its coalgebra structure is defined as 
\begin{align*}
&  \Delta(t_{ij}^{\pm}(u))=\sum_{k=1}^N
   t_{kj}^{\pm}(u\pm \frac{1}{4}\hbar c_2)\otimes
   t_{ik}^{\pm}(u\mp \frac{1}{4}\hbar c_1), \\             
&  \vep(T^{\pm}(u))=I~,~S({}^{t}T^{\pm}(u))=[{}^{t}T^{\pm}(u)]^{-1}, \\
&  \Delta(c)=c\otimes 1+1 \otimes c~,~\vep(c)=0~,~S(c)=-c,
\end{align*}
where $c_1=c\otimes 1$ and $c_2=1\otimes c$.
\end{defn}
Note that the subalgebra generated by
$\{ t_{ij}^k|1\leq i,j\leq N~,~k\in \nz_{\geq 0} \}$ is  
$Y_{\hbar}(\gtgl_N)$~ \cite{dr_1,dr_2} and the algebra 
$\cD Y_{\hbar}(\gtgl_N)$ is the quantum double of $Y_{\hbar}(\gtgl_N)$.
Let us define the pairing  $\langle \cdot, \cdot \rangle $ between
$T^{\pm}(u)$ as follows (cf. \cite{QISM}).
\[ \langle T^{+}(u), T^{-}(v) \rangle :=
   \sum_{i,j,k,l}\langle t^{+}_{ij}(u),t^{-}_{kl}(v) \rangle
   E_{ij}\otimes E_{kl} = R(u-v). \]
It seems that the following theorem is well-known to the specialists.
\begin{thm}
The pairing $\langle \cdot, \cdot \rangle $ gives the Hopf pairing.
\end{thm}  
The crucial point of the theorem is its non degeneracy. We could check
the non degeneracy for $N=2$ directly.
For the motivation of our choice, see Appendix \ref{quantum_affine}.
\subsection{Drinfel'd generators}
 
 We introduce the Drinfel'd generators
of $\cD Y_{\hbar}(\gtgl_N)$ exactly in the same way as in
\cite{ding}. 
\begin{thm}\label{ding_frenkel}
           $T^{\pm}(u)$ have the following unique decompositions:
\begin{align*}
T^{\pm}(u) &= \begin{pmatrix}
      1 & & & 0 \\
      f_{2,1}^{\pm}(u) & \ddots & & \\
            &\ddots &\ddots & \\
      f_{N,1}^{\pm}(u)& &f_{N,N-1}^{\pm}(u) & 1
              \end{pmatrix}
      \begin{pmatrix}
      k_{1}^{\pm}(u) & & & 0 \\
              & \ddots & &   \\
              & & \ddots & \\
      0 & & & k_{N}^{\pm}(u)
              \end{pmatrix} \\
          & \times \begin{pmatrix}
            1 & e_{1,2}^{\pm}(u)& & e_{1,N}^{\pm}(u) \\
              & \ddots &\ddots & \\
              & & \ddots & e_{N-1,N}^{\pm}(u) \\
            0 & & & 1
            \end{pmatrix}.
\end{align*}
\end{thm}
To prove this theorem, we have only to show that each component 
$f^{\pm}_{p,q}(u),k^{\pm}_{p}(u),e^{\pm}_{p,q}(u)$ is well-defined element
of $\cD Y_{\hbar}(\gtgl_N)[[u^{\mp 1}]]$. From the explicit formulae of these 
elements in terms of quantum minors, given in Appendix \ref{gauss}, it 
immediately follows since our algebra $\cD Y_{\hbar}(\gtgl_N)$ is 
$\hbar$-adically completed.  
Set
\begin{align*}
\hspace{1in} &X_{i}^-(u)=f_{i+1,i}^+(u_+)-f_{i+1,i}^-(u_-), \\
             &X_{i}^+(u)=e_{i,i+1}^+(u_-)-e_{i,i+1}^-(u_+).
\end{align*}
They satisfy the following commutation relations:
\begin{thm}\label{gl_N}
\begin{align*}
             & k_{i}^{\pm}(u)k_{j}^{\pm}(v)
              =k_{j}^{\pm}(v)k_{i}^{\pm}(u),~
               k_{i}^{+}(u)k_{i}^{-}(v)=k_{i}^{-}(v)k_{i}^{+}(u), \\
             & \frac{u_{\mp}-v_{\pm}}{u_{\mp}-v_{\pm}+\hbar}
               k_{i}^{\mp}(v)^{-1}k_{j}^{\pm}(u)=
               \frac{u_{\pm}-v_{\mp}}{u_{\pm}-v_{\mp}+\hbar} 
               k_{j}^{\pm}(u)k_{i}^{\mp}(v)^{-1}
               \quad i>j~, \\
             & \begin{cases}
               k_{i}^{\pm}(u)^{-1}X_{i}^{+}(v)k_{i}^{\pm}(u)=
               \dfrac{u_{\pm}-v+\hbar}{u_{\pm}-v}X_{i}^{+}(v), &  \\
               k_{i}^{\pm}(u)X_{i}^{-}(v)k_{i}^{\pm}(u)^{-1}=
               \dfrac{u_{\mp}-v+\hbar}{u_{\mp}-v}X_{i}^{-}(v), &                \end{cases} \\
             & \begin{cases}
               k_{i+1}^{\pm}(u)^{-1}X_{i}^{+}(v)k_{i+1}^{\pm}(u)=               \dfrac{u_{\pm}-v-\hbar}{u_{\pm}-v}X_{i}^{+}(v), & \\
               k_{i+1}^{\pm}(u)X_{i}^{-}(v)k_{i+1}^{\pm}(u)^{-1}=               \dfrac{u_{\mp}-v-\hbar}{u_{\mp}-v}X_{i}^{-}(v), &                \end{cases} \\
             & k_{j}^{\pm}(u)^{-1}X_{i}^{+}(v)k_{j}^{\pm}(u)
              =X_{i}^{+}(v),~
               k_{j}^{\pm}(u)X_{i}^{-}(v)k_{j}^{\pm}(u)^{-1}
              =X_{i}^{-}(v)
               \quad \text{otherwise}, \\
             & (u-v\mp \hbar)X_{i}^{\pm}(u)X_{i}^{\pm}(v)=
               (u-v\pm \hbar)X_{i}^{\pm}(v)X_{i}^{\pm}(u), \\
             & (u-v+\hbar)X_{i}^+(u)X_{i+1}^+(v)=
               (u-v)X_{i+1}^+(v)X_{i}^+(u), \\
             & (u-v)X_{i}^-(u)X_{i+1}^-(v)=
               (u-v+\hbar)X_{i+1}^-(v)X_{i}^-(u), \\
             & X_{i}^{\pm}(u_1)X_{i}^{\pm}(u_2)X_{j}^{\pm}(v)-
              2X_{i}^{\pm}(u_1)X_{j}^{\pm}(v)X_{i}^{\pm}(u_2)+
               X_{j}^{\pm}(v)X_{i}^{\pm}(u_1)X_{i}^{\pm}(u_2) \\
             & +\{ u_1 \leftrightarrow u_2 \}=0 \quad |i-j|=1, \\             & X_{i}^{\pm}(u)X_{j}^{\pm}(v)=
               X_{j}^{\pm}(v)X_{i}^{\pm}(u)\quad |i-j|>1, \\
             & [X_{i}^+(u),X_{j}^-(v)]=\hbar \delta_{ij}
               \left\{ \delta(u_{-}-v_{+})
                       k_{i+1}^+(u_{-})k_{i}^+(u_{-})^{-1}-
                       \delta(u_{+}-v_{-})
                       k_{i+1}^-(v_{-})k_{i}^-(v_{-})^{-1} \right\}.
\end{align*}
Here $\delta(u-v)=\sum_{k\in \nz} u^{-k-1}v^k$ is a delta function.
\end{thm}
One can prove the above theorem in exactly the same way 
as in \cite{ding} for the $U_{q}(\hgtgl_n)$ case.

\subsection{Two subalgebras}

 To decompose $\cD Y_{\hbar}(\gtgl_N)$
into two subalgebras 
$\cD Y_{\hbar}(\gtsl_N)$ and a Heisenberg subalgebra,
we introduce the following currents:
\begin{align*}
&  H_i^{\pm}(u)=k_{i+1}^{\pm}(u+\frac{1}{2}\hbar i)
              k_{i}^{\pm}(u+\frac{1}{2}\hbar i)^{-1}~,~
   K^{\pm}(u)=\prod_{i=1}^N 
   k_i^{\pm}\left( u+(i-\frac{N+1}{2})\hbar\right), \\
&  E_i(u)=\frac{1}{\hbar}X_i^{+}(u+\frac{1}{2}\hbar i)~,~
   F_i(u)=\frac{1}{\hbar}X_i^{-}(u+\frac{1}{2}\hbar i). 
\end{align*}
We define $\cD Y_{\hbar}(\gtsl_N)$ to be the subalgebra of
 $\cD Y_{\hbar}(\gtgl_N)$ generated by
$H_i^{\pm}(u)~,~E_i(u)~,~ F_i(u)~$ and $c$. A Heisenberg subalgebra 
of $ \cD Y_{\hbar}(\gtgl_N) $ generated by $K^{\pm}(u)$ 
commute with all of the elements of $ \cD Y_{\hbar}(\gtsl_N) $.
In fact we see that the formula $K^{\pm}(u)=\qdet.T^{\pm}(u)$ holds 
as a consequence of Theorem \ref{gauss_decomp}. (See Appendix \ref{gauss}
for the definition of $\qdet.T^{\pm}(u)$.) 
In terms of these generators, the above commutation relations can be 
rephrased as follows. Let $A=(a_{ij})$ be the Cartan matrix of the
Lie algebra $\gtsl_N$.
\begin{cor}
\begin{align*}\label{sl_N}
& [H_{i}^{\pm}(u),H_{j}^{\pm}(v)]=0, \\
&  (u_{\mp}-v_{\pm}+\hbar B_{ij})(u_{\pm}-v_{\mp}-\hbar B_{ij})
  H_{i}^{\pm}(u)H_{j}^{\mp}(v) \\
& =(u_{\mp}-v_{\pm}-\hbar B_{ij})(u_{\pm}-v_{\mp}+\hbar B_{ij})
  H_{j}^{\mp}(v)H_{i}^{\pm}(u), \\
& [K^{\pm}(u),K^{\pm}(v)]=0, \\
& f(u_{-}-v_{+})K^{+}(u)K^{-}(v)=K^{-}(v)K^{+}(u)f(u_{+}-v_{-}), \\
& \begin{cases}
  H_{i}^{\pm}(u)^{-1}E_{j}(v)H_{i}^{\pm}(u)=
  \dfrac{u_{\pm}-v-\hbar B_{ij}}{u_{\pm}-v+\hbar B_{ij}}E_{j}(v), & \\
  H_{i}^{\pm}(u)F_{j}(v)H_{i}^{\pm}(u)^{-1}=
  \dfrac{u_{\mp}-v-\hbar B_{ij}}{u_{\mp}-v+\hbar B_{ij}}F_{j}(v), &
  \end{cases} \\
& [K^{\sigma}(u),H_{i}^{\pm}(v)]=[K^{\sigma}(u),E_{i}(v)]
 =[K^{\sigma}(u),F_{i}(v)]=0 \quad \forall \sigma=\pm, \forall i \\
& (u-v-\hbar B_{ij})E_i(u)E_j(v)=(u-v+\hbar B_{ij})E_j(v)E_i(u), \\
& (u-v+\hbar B_{ij})F_i(u)F_j(v)=(u-v-\hbar B_{ij})F_j(v)F_i(u), \\
& \sum_{\sigma \in \gtS_m}
  [F_i(u_{\sigma(1)}),[F_i(u_{\sigma(2)})\cdots ,
  [F_i(u_{\sigma(m)}),F_j(v)]\cdots ]=0\quad i\neq j~,~m=1-a_{ij}, \\
& \sum_{\sigma \in \gtS_m}
  [E_i(u_{\sigma(1)}),[E_i(u_{\sigma(2)})\cdots ,
  [E_i(u_{\sigma(m)}),E_j(v)]\cdots ]=0\quad i\neq j~,~m=1-a_{ij}, \\
& [E_i(u),F_j(v)]=\frac{1}{\hbar}\delta_{ij}
  \left\{ \delta(u_{-}-v_{+})H_i^{+}(u_{-})-
          \delta(u_{+}-v_{-})H_i^{-}(v_{-})\right\}.
\end{align*}
Here we have set $B_{ij}=\frac{1}{2}(\alpha_i,\alpha_j)$ and
\[ f(u)=\prod_{j=1}^{N-1}\frac{u-j\hbar}{u+j\hbar}. \]
\end{cor}

To compare with the known results at $c=0$ \cite{tolstoy},
let us write down the commutation relations componentwise.
 The Fourier components of the generating series
$H_{i}^{\pm}(u),E_{i}(u),F_{i}(u)$ are of the following form:
\begin{align*}
& H_{i}^{+}(u)=1+\hbar\sum_{k\geq 0}h_{ik}u^{-k-1},~
  H_{i}^{-}(u)=1-\hbar\sum_{k<0}h_{ik}u^{-k-1}, \\
& E_{i}(u)=\sum_{k\in \nz}e_{ik}u^{-k-1},~
  F_{i}(u)=\sum_{k\in \nz}f_{ik}u^{-k-1}.
\end{align*}
For $c=0$, the commutation relations of $ \cD Y_{\hbar}(\gtsl_N)$
in terms of the above Fourier component look simple as follows:
\begin{cor}\label{comm_0}
\begin{align*}
& [h_{ik},h_{jl}]=0,~ [h_{i0},x_{jl}^{\pm}]=\pm2B_{ij}x_{jl}^{\pm},~
  [x_{ik}^+,x_{jl}^-]=\delta_{ij}h_{ik+l}, \\
& [h_{ik+1},x_{jl}^{\pm}]-[h_{ik},x_{jl+1}^{\pm}]=
  \pm \hbar B_{ij}[h_{ik},x_{jl}^{\pm}]_{+}, \\
& [x_{ik+1}^{\pm},x_{jl}^{\pm}]-[x_{ik}^{\pm},x_{jl+1}^{\pm}]=
  \pm \hbar B_{ij}[x_{ik}^{\pm},x_{jl}^{\pm}]_{+}, \\
& \sum_{\sigma\in \gtS_m}
  [x_{ik_{\sigma(1)}}^{\pm},[x_{ik_{\sigma(2)}}^{\pm},\cdots,
  [x_{ik_{\sigma(m)}}^{\pm},x_{jl}^{\pm}]\cdots ]=0 \quad
  i\neq j~,~m=1-a_{ij},
\end{align*}
for $k,l \in \nz$, where we set $x_{ik}^{+}=e_{ik}~,~x_{ik}^{-}=f_{ik}$ 
and $[x,y]_{+}=xy+yx$ for $x,y \in \cD Y_{\hbar}(\gtsl_N)$.
\end{cor}
These relations are the same as in \cite{tolstoy} with $\hbar=1$.
The set $\{ h_{ik}, x_{ik}^{\pm}|1\leq i \leq N-1, k \in \nz_{\geq 0} \}$ 
provides the Drinfel'd generators of $Y_{\hbar}(\gtsl_N)$ \cite{dr_2}.

Let us set
\[
E_{i}^{\pm}(u)=
\frac{1}{\hbar}e_{i,i+1}^{\pm}(u_{\mp}+\frac{1}{2}\hbar i),
\quad 
F_{i}^{\pm}(u)=
\frac{1}{\hbar}f_{i+1,i}^{\pm}(u_{\pm}+\frac{1}{2}\hbar i),
\]
so that $E_{i}^{\pm}(u)=E_{i}^+(u)-E_{i}^-(u),
         F_{i}^{\pm}(u)=F_{i}^+(u)-F_{i}^-(u)$. 
(See Theorem \ref{ding_frenkel} for the definition of 
 $e_{i,i+1}^{\pm}(u)$ and $f_{i+1,i}^{\pm}(u)$. ) 
Let us denote $\cD Y$ for $\cD Y_{\hbar}(\gtsl_N)$ and $\cD Y^{\pm}$ 
be the subalgebra of $\cD Y$ generated by $e_{i,k},f_{i,k}$ respectively.
 Set
\[ \cN^{\pm}=\sum_{i,k}x_{ik}^{\pm}\cD Y^{\pm}. \]
We get the partial results of the coproduct of the 
currents $E_{i}^{\pm}(u),F_{i}^{\pm}(u),H_{i}^{\pm}(u),K^{\pm}(u)$ which is
sufficient for our purpose.
\begin{lemma}\label{coproduct}
\begin{align*}
1)~ \Delta(E_{i}^{\pm}(u)) & \equiv 
     E_{i}^{\pm}(u)\otimes 1 + 
     H_{i}^{\pm}(u_{\mp})\otimes E_{i}^{\pm}(u\mp \frac{1}{2}\hbar c_1), \\
2)~ \Delta(F_{i}^{\pm}(u)) & \equiv
     1\otimes F_{i}^{\pm}(u) +
     F_{i}^{\pm}(u \pm \frac{1}{2}\hbar c_2)\otimes H_{i}^{\pm}(u_{\pm}), \\
3)~ \Delta(H_{i}^{\pm}(u)) & \equiv
     H_{i}^{\pm}(u\pm \frac{1}{4}\hbar c_2)\otimes
     H_{i}^{\pm}(u\mp \frac{1}{4}\hbar c_1), \\
  & \mod (\cN^{-} \cD Y \otimes \cD Y \cN^{+})\cap
         (\cD Y \cN^{-} \otimes \cN^{+} \cD Y), \quad \text{and} \\
4)~ \Delta(K^{\pm}(u)) & =
     K^{\pm}(u\pm \frac{1}{4}\hbar c_2)\otimes 
     K^{\pm}(u\mp \frac{1}{4}\hbar c_1).
\end{align*}
\end{lemma} 
The last formula follows from the fact that the equation 
$K^{\pm}(u)=\qdet.T^{\pm}(u)$ holds.
We remark that these formula for $Y_{\hbar}(\gtsl_N)$ are obtained in
\cite{chari2}. The exact formulae in the case of $\cD Y_{\hbar}(\gtgl_2)$
are given in \cite{kohno}. For more information on the coproduct formulae,
see Appendix \ref{fusion_T}.

\subsection{Quantum current}\label{current}
Here we give a remark of the Definition \ref{defn_DY(gl_N)} and introduce
the so-called quantum current \cite{resh}.

Let $\overline{R}(u)$ be the Yang's $R$-matrix normalized as 
(\ref{Yang's R}) and $R(u)=f_{N}(u)\overline{R}(u)$ 
with some scalar function $f_{N}(u)$.
We remark that even if we change the normalization of $R$-matrix in 
Definition \ref{defn_DY(gl_N)} to $R(u)$ defined here, 
the commutation relations given by Corollary \ref{sl_N} never change 
except for the relations between $K^{\pm}(u)$. It changes as follows.
\[ f(u_{-}-v_{+})K^{+}(u)K^{-}(v)=K^{-}(v)K^{+}(u)f(u_{+}-v_{-}), \]
where
\[ f(u)=\{\prod_{i,j=1}^N f_N(u+(i-j)\hbar)\}
        \prod_{k=1}^{N-1} \frac{u-k\hbar}{u+k\hbar}. \]
If we set
\[ f_N(u)=\frac{\Gamma(\frac{u}{N\hbar})\Gamma(1+\frac{u}{N\hbar})}
               {\Gamma(\frac{1}{N}+\frac{u}{N\hbar})
                \Gamma(1-\frac{1}{N}+\frac{u}{N\hbar})} \] 
 where $\Gamma(u)$ is the Euler's Gamma function, then $f(u)=1$.   
In the rest of this subsection, we fix the function $f_N(u)$ as above.
Let us define the Quantum current $T(u)$ as
\[ T(u)=T^{+}(u_{-})T^{-}(u_{+})^{-1}. \]
They enjoy the following commutation relations.
\begin{lemma}
\begin{align*}
& R(u-v)\overset{1}{T}(u)R_{21}(v-u-\hbar c)\overset{2}{T}(v)=
  \overset{2}{T}(v)R(u-v-\hbar c)\overset{1}{T}(u)R_{21}(v-u), \\
& R(u_{-}-v_{+})\overset{1}{T}(u)R_{21}(v_{-}-u_{+})
  \overset{2}{T^{\pm}}(v_{\pm})=
  \rho_N(u_{\mp}-v_{\pm})\overset{2}{T^{\pm}}(v_{\pm})\overset{1}{T}(u),
\end{align*}
where $\rho_N(u)=f_N(u)f_N(-u)$.
\end{lemma}
Since $T(u)$ can be regarded as $N\times N$ matrix, 
we can define the current $l(u)$ by
\[ l(u)=\text{{\rm tr}.}T(u). \]
At the critical level ($c=-N$), one can show that $l(u)$ commutes with
$T^{\pm}(u)$ so that  $l(u)$ provides the Yangian deformed Gelfand Dickii
algebra \cite{Vir}.
\begin{rem} Everything given in this section makes sense except for this 
subsection. Since the function $f_N(u)$ chosen here can not be regarded
as a formal series in $\hbar$, the formulae given after the specific 
choice of $f_N(u)$ must be considered only as an asymptotics.  
\end{rem}  
\section{Representation theory of $\cD Y_{\hbar}(\gtg)$ }\hspace{1 in} \\
 Unfortunately, we have no general theorem about the representation theory of
$\cD Y_{\hbar}(\gtg)$ at the moment due to the lack of triangular decomposition
and the grading operator $d$. Nevertheless, we expect that the representation 
theory of $\cD Y_{\hbar}(\gtg)$ can be established just as in the case of 
quantum affine algebra \cite{guide,nankai}.

 In this section we present examples of both finite and infinite
dimensional representations of $\cD Y_{\hbar}(\gtg)$.
\subsection{Finite dimensional representations}\label{FDR}
At $c=0$, the Heisenberg subalgebra becomes central in 
$\cD Y_{\hbar}(\gtgl_N)$. On behalf of Schur's Lemma, it is sufficient for
investigating the irreducible finite dimensional representations to
consider $\cD Y_{\hbar}(\gtsl_N)$ case. From Corollary \ref{comm_0}, 
we expect that most of the finite dimensional $Y_{\hbar}(\gtsl_N)$-module 
can be extended to $\cD Y_{\hbar}(\gtsl_N)$-module. 

Let {\bf d}$=\{d_{i,k}\}_{1\leq i\leq N-1, k\in \nz}$ be a subset of $\cA$.
A $\cD Y_{\hbar}(\gtsl_N)$-module $V$ is called 
{\it pseudo highest weight module with pseudo highest weight {\bf d}} if there
is an unique, up to scalar multiple, non-zero vector $v\in V$ such that $V$
is generated by $v$ and
\[ h_{i,k}.v=d_{i,k}v, \quad e_{i,k}.v=0, \qquad 1\leq \forall i \leq N-1,~
                                          \forall k \in \nz. \]
Here we borrow this terminology from \cite{guide}.
Let us denote $V(${\bf d}$)$ for such $V$.
\begin{conj} \hspace{1 in} \\
\begin{enumerate}
\renewcommand{\labelenumi}{(\roman{enumi})}
\item Let $V$ be an irreducible finite dimensional 
      $Y_{\hbar}(\gtsl_N)$-module whose constant term of the Drinfel'd 
      polynomials are invertible. Then $V$ can be lift up to an irreducible
      finite dimensional $\cD Y_{\hbar}(\gtsl_N)$-module.
\item The irreducible $\cD Y_{\hbar}(\gtsl_N)$-module $V(${\bf d}$)$ of pseudo
      highest weight {\bf d} is finite dimensional iff there exist monic 
      polynomials $P_{i}(v) \in \cA[v]~1\leq i\leq N-1$ such that
\[ 1-\hbar\sum_{k<0}d_{ik}v^{-k-1}=
   \frac{P_{i}(v+\frac{1}{2}(\alpha_i,\alpha_i)\hbar)}{P_{i}(v)}=
   1+\hbar\sum_{k\geq 0}d_{ik}v^{-k-1}, \]
      in the sense that the left-hand side and the right-hand side are the 
      Laurent expansion of the middle term
      about $0$ and $\infty$ respectively.
\end{enumerate}
\end{conj}
The above monic polynomials $P_i$ are called Drinfel'd polynomials.
Next we show some examples which support this conjecture. 

\begin{ex}{\bf $\gtsl_2$ case} 

Here we omit writing the subscript $1$ for simplicity.
Let $W_m=\oplus_{j=0}^m \cA w_j$ be the spin $\frac{m}{2}$ representation of
$\gtsl_2$ and set
\[ W_m(u)=W_m \otimes_{\cA} \cA((u^{-1})), \]
where $u$ is thought to be either an indeterminate or an invertible 
element of $\cA$.
It is known by \cite{chari1} that we can define the  
$Y_{\hbar}(\gtsl_2)$-module structure on $W_m(u)$. It immediately follows
that their action of $Y_{\hbar}(\gtsl_2)$ can be extended to 
$\cD Y_{\hbar}(\gtsl_2)$.

\begin{lemma}The action $\cD Y_{\hbar}(\gtsl_2)$ is given by
\begin{align*}
(1)\qquad &  e_k.w_i=\left\{ u+(\frac{1}{2}m-i+\frac{1}{2})\hbar \right\}^k
                    (m-i+1)w_{i-1}, \\
(2)\qquad &  f_k.w_i=\left\{ u+(\frac{1}{2}m-i-\frac{1}{2})\hbar \right\}^k
                    (i+1)w_{i+1}, \\
(3)\qquad &  h_k.w_i=\left[ 
                     \left\{ u+(\frac{1}{2}m-i-\frac{1}{2})\hbar \right\}^k
                     (i+1)(m-i) \right. \\
          & \hspace{0.5 in} \left. -
                     \left\{ u+(\frac{1}{2}m-i+\frac{1}{2})\hbar \right\}^k
                     i(m-i+1) \right] w_i, 
\end{align*}
where we set $w_{-1}=w_{m+1}=0$.
\end{lemma}
As a Consequence, we obtain the following.
\begin{cor} \hspace{1 in} \\
\begin{enumerate}
\renewcommand{\labelenumi}{(\roman{enumi})}
\item $W_{m}(u)$ is a pseudo highest weight module with pseudo highest weight
      {\bf d}=$\{ d_{k}\}$ given by
\[    d_k=m(u+\frac{m-1}{2}\hbar )^k . \]
\item The Drinfel'd polynomial $P$ associated to $W_{m}(u)$ is given by
\[ P(v)=\{v-u-\frac{m-1}{2}\hbar \}\{v-u-\frac{m-3}{2}\hbar \}\cdots
        \{v-u+\frac{m-1}{2}\hbar \}. \]
\end{enumerate}
\end{cor}
\end{ex}

\begin{ex}{\bf $\gtsl_N$ case~(vector representation)}\label{DY(sl_N)}

Let $u$ be either an indeterminate or an invertible element of $\cA$.
Set
\[ V_{u}=V\otimes_{\cA} \cA ((u^{-1})),\quad V=\oplus_{j=0}^{N-1}\cA w_j. \]
We can extend the action of $Y_{\hbar}(\gtsl_N)$ to $\cD Y_{\hbar}(\gtsl_N)$
as follows. (See Example \ref{Y(sl_N)}.)
\begin{lemma} The action of $\cD Y(\gtsl_N)$ is given by
\begin{align*}
(1)\qquad & e_{i,k}.w_{i}=(u-\frac{N-1-i}{2}\hbar)^k.w_{i-1},\qquad
            e_{i,k}.w_{j}=0 \quad j\neq i, \\
(2)\qquad & f_{i,k}.w_{i-1}=(u-\frac{N-1-i}{2}\hbar)^k.w_{i},\qquad
            f_{i,k}.w_{j}=0 \quad j\neq i-1, \\
(3)\qquad & h_{i,k}.w_{i-1}=(u-\frac{N-1-i}{2}\hbar)^k.w_{i-1}, \qquad
            h_{i,k}.w_{i}=-(u-\frac{N-1-i}{2}\hbar)^k.w_{i}, \\
          & h_{i,k}.w_{j}=0 \quad j\neq i,i-1, 
\end{align*}
\end{lemma}
Hence we have the following.
\begin{cor} \hspace{1 in} \\
\begin{enumerate}
\renewcommand{\labelenumi}{(\roman{enumi})}
\item $V_{u}$ is a pseudo highest weight module with highest weight
      {\bf d}=$\{d_{ik}\}$ given by
\[   d_{1k}=(u-\frac{N-2}{2}\hbar)^k,\qquad d_{ik}=0 \quad i\neq 1. \]
\item The Drinfel'd polynomials $P_i$ associated to $V_u$ are given by
\[   P_{1}(v)=v-(u-\frac{N-2}{2}\hbar),\qquad P_{i}(v)=1 \quad i\neq 1. \]
\end{enumerate}
\end{cor}
\end{ex}

The next example is the generalization of the above example.

\begin{ex}{\bf $\gtsl_N$ case} Let $\gtg$ be a Lie algebra of type $A_{N-1}$.
In Example \ref{generic_repr}, we define irreducible finite dimensional
$Y_{\hbar}(\gtg)$-modules $V_a(m\Lambda_i)$ for $1\leq i\leq N-1$.
Here we give a sketch of proof that we can extend its action to 
$\cD Y_{\hbar}(\gtg)$. To see this
\begin{enumerate}
\renewcommand{\labelenumi}{(\roman{enumi})}
\item Define $\cD Y_{\hbar}(\gtg)$-module structure for $m=1,\forall i$. 
( Calculate the action on each weight vector explicitly, 
then it turns out that the invertibilty in the conjecture is essential.) 
\item Using the following embedding of $Y_{\hbar}(\gtg)$-module
\[ V_a((m+1)\Lambda_i)  \hookrightarrow
   V_{a-\frac{1}{2}\hbar}(m\Lambda_i)\otimes V_{a+\frac{m}{2}\hbar}(\Lambda_i),
\]
prove that we can define $\cD Y_{\hbar}(\gtg)$-module structure
for $\forall m, \forall i$ by induction on $m$. 
\end{enumerate}
The details are left to the reader as an excersize. 

 Note that the Drinfel'd polynomials of $V_a(m\Lambda_i)$ are exactly 
the same as in Example \ref{generic_repr}.
\end{ex}

\subsection{Bosonization of level 1 module}\label{bosonization}

Here we construct level $1$ $ \cD Y_{\hbar}(\gtg)$-module and 
vertex operators for $\gtg=\gtgl_N,\gtsl_N$ directly in terms of bosons. 

Let $\gth=\oplus_{i=1}^N \nc \vep_i $ be a Cartan subalgebra of
$\gtgl_N$, $\overline{Q}=\oplus_{i=1}^{N-1} \nz \alpha_i~
(\alpha_i=\vep_i-\vep_{i+1})$ be the root lattice of $\gtsl_N$,  
$\overline{\Lambda}_i=\Lambda_i-\Lambda_0$ be the classical part of the
$i$-th fundamental weight and $(\cdot, \cdot)$ be the standard bilinear
form defined by $(\vep_i,\vep_j)=\delta_{ij}$.
Let us introduce bosons 
$\{ a_{i,k}|1\leq i \leq N,~k\in \nz \setminus \{ 0\} \}$ 
satisfying:
\[ [a_{i,k},a_{j,l}]=k\delta_{i,j}\delta_{k+l,0}. \]

\noindent{\bf 1. $\gtgl_N$ case} \hspace{1 in} \\

Set
\[ \cF_{i,s}:=\cA [a_{j,-k}(1\leq j\leq N,~k\in \nz_{>0})]\otimes
              \cA [\overline{Q}]
              e^{\overline{\Lambda}_i+\frac{s}{N}(\sum_{j=1}^N\vep_j)}
   \quad (0\leq i\leq N-1), \]
where $s$ is a complex parameter and $\cA[\overline{Q}]$ is 
the group algebra of $\overline{Q}$ over $\cA$. 
On this space, we define the action of the operators 
$a_{j,k},\partial_{\vep_j},e^{\vep_j}$~$(1\leq j \leq N)$ by
\begin{align*}
a_{j,k}\cdot f\otimes e^{\beta}& =\begin{cases}
                      a_{j,k}f \otimes e^{\beta}          & k< 0  \\
            \text{$[a_{j,k},f]$}  \otimes e^{\beta} \quad & k> 0
                                  \end{cases}, \\
\partial_{\vep_j}\cdot f\otimes e^{\beta}& 
        =(\vep_j,\beta)f\otimes e^{\beta}, 
         \qquad \text{for}~f\otimes e^{\beta}\in \cF_{i,s} \\
e^{\vep_j}\cdot f\otimes e^{\beta} & 
=f\otimes e^{\vep_j+\beta}.
\end{align*}

\begin{thm} \label{level 1 gl_N}
The following assignment defines a 
$ \cD Y_{\hbar}(\gtgl_N)$-module structure on $\cF_{i,s}$.
\begin{align*}
k_j^+(u)&\mapsto \exp \left[ -\sum_{k>0}
          \frac{a_{j,k}}{k}\left \{ (u+\frac{1}{2}\hbar)^{-k}
                                  -(u-\frac{1}{2}\hbar)^{-k}
          \right \} \right]
          \left( \frac{u-\frac{1}{2}\hbar}{u+\frac{1}{2}\hbar} 
          \right)^{\partial_{\vep_j}}, \\
k_j^-(u)&\mapsto \exp \left[ \sum_{k>0,r<j}\frac{a_{r,-k}}{k}
          \left \{ u^k-(u-\hbar)^k \right \}+
                      \sum_{k>0,r>j}\frac{a_{r,-k}}{k}
          \left \{ (u+\hbar)^k-u^k \right \} \right], \\
\frac{1}{\hbar}X_j^+(u)& \mapsto \exp \left[ -\sum_{k>0}
          \frac{a_{j,-k}}{k}(u-\frac{3}{4}\hbar)^k+
          \sum_{k>0}\frac{a_{j+1,-k}}{k}(u+\frac{1}{4}\hbar)^k
          \right] \\
               & \times \exp \left[ \sum_{k>0}
          \frac{a_{j,k}-a_{j+1,k}}{k}(u+\frac{1}{4}\hbar)^{-k}
          \right] e^{\alpha_{j}}
          \left[ (-1)^{j-1}(u+\frac{1}{4}\hbar)
          \right]^{\partial_{\alpha_j}}, \\
\frac{1}{\hbar}X_j^-(u)& \mapsto \exp \left[ \sum_{k>0}
          \frac{a_{j,-k}}{k}(u-\frac{1}{4}\hbar)^k-
          \sum_{k>0}\frac{a_{j+1,-k}}{k}(u+\frac{3}{4}\hbar)^k
          \right] \\
               & \times \exp \left[ \sum_{k>0}
          \frac{-a_{j,k}+a_{j+1,k}}{k}(u-\frac{1}{4}\hbar)^{-k}
          \right] e^{-\alpha_{j}}
          \left[ (-1)^{j-1}(u-\frac{1}{4}\hbar)
          \right]^{-\partial_{\alpha_j}}.
\end{align*}
where we set $\partial_{\alpha_{j}}=
              \partial_{\vep_{j}}-\partial_{\vep_{j+1}}$.
\end{thm}

 Next we present the bosonization of type $I$ and type $II$ 
vertex operators. For this purpose, 
let us consider the evaluation module. Set
\[ V_u=V\otimes_{\cA}\cA((u^{-1})),\quad
   V=\oplus_{j=0}^{N-1} \cA w_j. \]
We define the $ \cD Y_{\hbar}(\gtgl_N)$-module structure on $V_u$
as follows:
\begin{align*}
 k^{\pm}_{i+1}(v).w_i=f^{\pm}(v-u)
 \frac{v-u+(\frac{N-3}{2}-i)\hbar}{v-u+(\frac{N-1}{2}-i)\hbar}w_i,& \quad
 k^{\pm}_j(v).w_i=f^{\pm}(v-u)w_i \quad \text{otherwise},  \\
 X^{+}_i(v).w_i=\hbar \delta(v-u+(\frac{N-1}{2}-i)\hbar)w_{i-1},& \quad
 X^{+}_j(v).w_i=0 \quad \text{otherwise}, \\
 X^{-}_i(v).w_{i-1}=\hbar \delta(v-u+(\frac{N-1}{2}-i)\hbar)w_{i},& \quad
 X^{-}_j(v).w_i=0 \quad \text{otherwise}, 
\end{align*}
where we set 
\[ f^{+}(u)=1,\quad f^{-}(u)=\frac{u-\frac{N-1}{2}\hbar}{u+\frac{N-3}{2}\hbar}.
\]
We remark that the restriction of the action of $\cD Y_{\hbar}(\gtgl_N)$ 
on the above $V_u$ to that of $\cD Y_{\hbar}(\gtsl_N)$ gives $V_u$ in
Example \ref{DY(sl_N)} exactly. 
\begin{defn}Vertex operators are intertwiners of the following form:
\begin{enumerate}
\renewcommand{\labelenumi}{(\roman{enumi})}
\item (type $I$) \hspace{1 in} 
      $\displaystyle{ \Phi^{(i,i+1)}(u):\cF_{i+1,s}\longrightarrow
                      \cF_{i,s-1}\otimes V_{u}, }$ 
\item (type $II$) \hspace{0.95 in}
      $\displaystyle{ \Psi^{(i,i+1)}(u):\cF_{i+1,s}\longrightarrow
                      V_{u}\otimes \cF_{i,s-1}. }$
\end{enumerate}
Here the indices are considered modulo $N$.
\end{defn}
Set 
\[ \Phi^{(i,i+1)}(u)=\sum_{j=0}^{N-1} \Phi^{(i.i+1)}_j(u)\otimes w_{j},
\qquad
   \Psi^{(i,i+1)}(u)=\sum_{j=0}^{N-1} w_{j}\otimes \Psi^{(i.i+1)}_j(u),
\]
We normalize them as 
\begin{enumerate}
\renewcommand{\labelenumi}{(\roman{enumi})}
\item $\displaystyle{ 
       \langle \Lambda_{i},s-1|\Phi^{(i,i+1)}_{i}(u)|\Lambda_{i+1},s \rangle=1,
       } $
\item $\displaystyle{
       \langle \Lambda_{i},s-1|\Psi^{(i,i+1)}_{i}(u)|\Lambda_{i+1},s \rangle=1,
       } $
\end{enumerate}
where we set 
$|\Lambda_i,s\rangle=1\otimes
e^{\overline{\Lambda}_i+\frac{s}{N}(\sum_{j=1}^N \vep_j)}$.
We mean by $\langle \Lambda_{i},s-1|
     \Phi^{(i,i+1)}_{i}(u)|\Lambda_{i+1},s \rangle$ the 
coefficient of $|\Lambda_{i},s-1 \rangle$ of the element
$\Phi^{(i,i+1)}_{i}(u)|\Lambda_{i+1},s \rangle$, and similarly for
$\Psi^{(i,i+1)}_{i}(u)$. With the above normalization our vertex 
operators uniquely exist. By using Lemma \ref{coproduct}, 
we obtain the bosonization formula of these vertex operators as follows.
\begin{thm}[Bosonization of vertex operators] \label{VOp_gl_N}
           For $0 \leq i \leq N-1$,
\begin{align*}               
\Phi^{(i,i+1)}_{N-1}(u)&= \exp \left[ \sum_{k>0}
            \frac{a_{N,-k}}{k}\left(
            u+(\frac{N}{2}+\frac{1}{4} )\hbar \right)^k \right] \\
            & \times      \exp \left[ \sum_{k>0;1\leq j<N}
            \frac{a_{j,k}}{k} \left(
            u-(\frac{N}{2}-\frac{1}{4}-j )\hbar \right)^{-k} 
                            \right]  \\
                     &  \times e^{-\vep_N}
            \left[(-1)^{N-1}\left(
            u+(\frac{N}{2}-\frac{3}{4})\hbar \right)
            \right]^{\partial_{\overline{\Lambda}_{N-1}}
                    + \frac{N-i-1}{N}}
            (-1)^{\frac{1}{2}(N-i-1)(N+i-2)}, \\
\Phi^{(i,i+1)}_{k-1}(u)&= [\Phi^{(i,i+1)}_{k}(u), f_{k,0}], \\
\Psi^{(i,i+1)}_{0}(u)&= \exp \left[ \sum_{k>0}
            \frac{a_{1,-k}}{k}\left(
            u-(\frac{N}{2}-\frac{3}{4})\hbar \right)^k \right] \\
            & \times    \exp \left[ \sum_{k>0;1<j \leq N}
            \frac{a_{j,k}}{k} \left(
            u-(\frac{N}{2}+\frac{1}{4}-j)\hbar \right)^{-k} 
                            \right] \\
                     &  \times e^{-\vep_1}
            \left[ -\left(
            u-(\frac{N}{2}-\frac{7}{4})\hbar \right)
            \right]^{-\partial_{\overline{\Lambda}_1}
                    + \frac{N-i-1}{N}}
            (-1)^{\frac{1}{2}i(i+1)}, \\
\Psi^{(i,i+1)}_{k}(u)&= [\Psi^{(i,i+1)}_{k-1}(u), e_{k,0}].
\end{align*}
\end{thm}

\noindent{\bf 2. $\gtsl_N$ case} \hspace{1 in} \\

Here we keep the same notation as in $\gtgl_N$ case unless otherwise
stated. Set
\[ \cF_{i}:=\cA [a_{j,-k}(1\leq j\leq N-1,~k\in \nz_{>0})]\otimes
              \cA [\overline{Q}]
              e^{\overline{\Lambda}_i}
   \quad (0\leq i\leq N-1). \]
As in the previous subsection, we define the action of the operators
$a_{j,k},\partial_{\alpha_j},e^{\alpha_j}$~$(1\leq j \leq N-1)$ on $\cF_i$.
\begin{thm} \label{level 1 sl_N}
The following assignment defines a 
$ \cD Y_{\hbar}(\gtsl_N)$-module structure on $\cF_{i}$.
\begin{align*}
H_j^+(u)&\mapsto  \exp \left[ -\sum_{k>0}\frac{a_{j,k}}{k}
           \left\{ (u+\frac{1}{2}\hbar)^{-k}
                  -(u-\frac{1}{2}\hbar)^{-k}\right\} \right]          
           \left( \frac{u-\frac{1}{2}\hbar}
                       {u+\frac{1}{2}\hbar} \right)^{-\partial_{\alpha_j}}, \\
H_j^-(u)&\mapsto  \exp \left[ -\sum_{k>0}\frac{a_{j,-k}}{k}
           \left\{ (u+\hbar)^k-(u-\hbar)^k \right\} \right. \\
        & \hspace{0.5 in} \left.
           + \sum_{k>0}\frac{a_{j+1,-k}+a_{j-1,-k}}{k}
           \left\{ (u+\frac{1}{2}\hbar)^k-(u-\frac{1}{2}\hbar)^k \right\}
           \right], \\
E_j(u)&\mapsto    \exp \left[ \sum_{k>0}\frac{a_{j,-k}}{k}
           \left\{ (u+\frac{1}{4}\hbar)^k+(u-\frac{3}{4}\hbar)^k \right\}
           - \sum_{k>0}\frac{a_{j+1,-k}+a_{j-1,-k}}{k}
           (u-\frac{1}{4}\hbar)^k \right] \\
      &    \times \exp \left[ -\sum_{k>0}\frac{a_{j,k}}{k}
           (u+\frac{1}{4}\hbar)^{-k} \right] e^{\alpha_j}
           \left[ (-1)^{j-1}(u+\frac{1}{4}\hbar) 
           \right]^{\partial_{\alpha_j}}, \\
F_j(u)&\mapsto    \exp \left[ -\sum_{k>0}\frac{a_{j,-k}}{k}
           \left\{ (u+\frac{3}{4}\hbar)^k+(u-\frac{1}{4}\hbar)^k \right\}
           +\sum_{k>0}\frac{a_{j+1,-k}+a_{j-1,-k}}{k}
           (u+\frac{1}{4}\hbar)^k \right] \\
      &    \times \exp \left[ \sum_{k>0}\frac{a_{j,k}}{k}
           (u-\frac{1}{4}\hbar)^{-k} \right] e^{-\alpha_j}
           \left[ (-1)^{j-1}(u-\frac{1}{4}\hbar)
           \right]^{-\partial_{\alpha_j}}.
\end{align*}
\end{thm}

Before investigating the vertex operators, we shall give some remarks here.
Every field in Theorem \ref{level 1 gl_N},\ref{VOp_gl_N},\ref{level 1 sl_N}
make sense as a formal series in $\hbar$ if we use the binomial expansion
\[ (u+a\hbar)^k=\sum_{j\geq 0} \begin{pmatrix}
                               k \\ j 
                               \end{pmatrix} (a\hbar)^ju^{k-j}
   \qquad a \in \cA, k\in \nz.                                   \]
Now one can prove these Theorem by some routine calculations. Notice that
because of the artificial choice of the action of the Heisenberg
subalgebra, the bosonization of the vertex operators in the case of 
$\cD Y_{\hbar}(\gtgl_N)$ has such nice expression. For the 
$\cD Y_{\hbar}(\gtsl_N)$ case, as we will see soon, 
we have some subtle problem to bosonize the vertex operators.

To introduce the vertex operators of type $I$ and type $II$,
let us fix the evaluation module $V_u$ given in Example \ref{DY(sl_N)}.
\begin{defn}Vertex operators are intertwiners of the following form:
\begin{enumerate}
\renewcommand{\labelenumi}{(\roman{enumi})}
\item (type $I$) \hspace{1 in} 
      $\displaystyle{ \Phi^{(i,i+1)}(u):\cF_{i+1}\longrightarrow
                      \cF_{i}\otimes V_{u}, }$ 
\item (type $II$) \hspace{0.95 in}
      $\displaystyle{ \Psi^{(i,i+1)}(u):\cF_{i+1}\longrightarrow
                      V_{u}\otimes \cF_{i}. }$
\end{enumerate}
Here the indices are considered modulo $N$.
\end{defn}
Set 
\[ \Phi^{(i,i+1)}(u)=\sum_{j=0}^{N-1} \Phi^{(i.i+1)}_j(u)\otimes w_{j},
\qquad
   \Psi^{(i,i+1)}(u)=\sum_{j=0}^{N-1} w_{j}\otimes \Psi^{(i.i+1)}_j(u),
\]
We normalize them as 
\begin{enumerate}
\renewcommand{\labelenumi}{(\roman{enumi})}
\item $\displaystyle{ 
       \langle \Lambda_{i}|\Phi^{(i,i+1)}_{i}(u)|\Lambda_{i+1}\rangle=1,
       } $
\item $\displaystyle{
       \langle \Lambda_{i}|\Psi^{(i,i+1)}_{i}(u)|\Lambda_{i+1}\rangle=1,
       } $
\end{enumerate}
Just as in the case of $\cD Y_{\hbar}(\gtgl_N)$, 
our vertex operators uniquely exist and the bosonization 
formulae are as follows.
\begin{thm}[Bosonization of vertex operators] \label{VOp_sl_N}
           For $0 \leq i \leq N-1$,
\begin{align*}               
\Phi^{(i,i+1)}_{N-1}(u)   &=\lim_{n\rightarrow \infty}
                            \Phi^{(i,i+1)}_{N-1}(u)_n, \\
\Phi^{(i,i+1)}_{N-1}(u)_n &= \exp \left[ \sum_{k>0}
                            \frac{a_{N-1,-k}}{k}\left(
                            u+\frac{3}{4} \hbar \right)^k \right]
                             \exp \left[ -\sum_{k>0;1\leq j<N}
                            \frac{a_{j,k}}{2k}
                            f_{j,k}^{(I)}(u-\frac{1}{4}\hbar)_n \right] \\
                          &\times e^{\overline{\Lambda}_{N-1}}
\left\{\prod_{j=1}^{N-1}[g_j^{(I)}(u-\frac{1}{4}\hbar)_n]^{\partial_{\alpha_j}}
\right\}  
[(-1)^{\frac{N-1}{2}}(N\hbar)^{-\frac{N-2}{2}}
]^{\partial_{\overline{\Lambda}_{N-1}}-\frac{i+1}{N}} \\
                          &\times g_{i+1}^{(I)}(u-\frac{1}{4}\hbar)_n^{-1}
[u-(\frac{N-i}{2}-\frac{3}{4})\hbar](-1)^{\frac{1}{2}(N-i)(N+i-1)}, \\
\Phi^{(i,i+1)}_{k-1}(u)   &= [\Phi^{(i,i+1)}_{k}(u), f_{k,0}], \\
\Psi^{(i,i+1)}_{0}(u)     &=\lim_{n\rightarrow \infty}
                            \Psi^{(i,i+1)}_{0}(u)_n, \\
\Psi^{(i,i+1)}_{0}(u)_n   &= \exp \left[ -\sum_{k>0}
                            \frac{a_{1,-k}}{k}\left( u-
                           (\frac{N}{2}-\frac{1}{4})\hbar \right)^k \right]
                            \exp \left[ \sum_{k>0;1\leq j<N}
                            \frac{a_{j,k}}{2k}
                            f_{N-j,k}^{(II)}(u+\frac{1}{4}\hbar)_n \right] \\
                          &\times e^{-\overline{\Lambda}_{1}}
\left\{\prod_{j=1}^{N-1}[g_{N-j}^{(II)}
      (u+\frac{1}{4}\hbar)_n]^{-\partial_{\alpha_j}}\right\}
[(-1)^{-\frac{1}{2}}(N\hbar)^{-\frac{N-2}{2}}]^{
-\partial_{\overline{\Lambda}_{1}}+\frac{N-i-1}{N}} \\
                          &\times g_{N-i-1}^{(II)}(u+\frac{1}{4}\hbar)_n
                           (-1)^{\frac{1}{2}i(i+1)}, \\ 
\Psi^{(i,i+1)}_{k}(u)&= [\Psi^{(i,i+1)}_{k-1}(u), e_{k,0}].
\end{align*}
Here the functions $f^{\ast}_{j,k}(u)_n,g^{\ast}_{j}(u)_n~(\ast=(I),(II))$
are defined as follows. 
\begin{align*}
&  f^{\ast}_{j,k}(u)_n=\sum_{l=0}^{j-1}f_{k}^{\ast}(u+\frac{j-1-2l}{2}\hbar)_n,
   \qquad 1\leq j< N, \\
&  g^{\ast}_{j}(u)_n=\begin{cases}
   \qquad  \displaystyle{1} & j=0 \\
   \displaystyle{\left[ \prod_{l=0}^{j-1}g^{\ast}
          (u+\frac{j-1-2l}{2}\hbar)_n \right]^{\frac{1}{2}}}
      & j>0 \end{cases}\quad , \\
& f_k^{(I)}(u)_n=(u-\frac{N-2}{2}\hbar)^{-k}+
  \sum_{l=0}^{n-1}\left\{(u+\frac{N}{2}\hbar+N\hbar l)^{-k}
                        -(u+\frac{N}{2}\hbar+(Nl+1)\hbar )^{-k} \right. \\
& \hspace{2.2 in} \left.
                        +(u-\frac{3N}{2}\hbar-(Nl-1)\hbar)^{-k}
                        -(u-\frac{N}{2}\hbar-N\hbar l)^{-k}\right\}, \\
& f_k^{(II)}(u)_n=u^{-k}+
  \sum_{l=0}^{n-1}\left\{(u+N\hbar+N\hbar l)^{-k}
                        -(u+\hbar+N\hbar l)^{-k} \right. \\
& \hspace{1.5 in} \left.
                        +(u-(N-1)\hbar-N\hbar l)^{-k}
                        -(u-N\hbar-N\hbar l)^{-k}\right\}, \\
& g^{(I)}(u)_n=(u-\frac{N-2}{2}\hbar)e^{\frac{N-2}{N}\gamma}
  \frac{(u+\frac{N}{2}\hbar)(u-\frac{3N}{2}\hbar+\hbar)}
       {(u-\frac{N}{2}\hbar)(u+\frac{N}{2}\hbar+\hbar)} \\
& \hspace{0.5 in}\times  \prod_{l=1}^{n-1}\left[ 
  \frac{(u+\frac{N}{2}\hbar+N\hbar l)(u-\frac{3N}{2}\hbar-(Nl-1)\hbar)}
       {(u-\frac{N}{2}\hbar-N\hbar l)(u+\frac{N}{2}\hbar+(Nl+1)\hbar )}
  \right]e^{-\frac{N-2}{Nl}}, \\
& g^{(II)}(u)_n=ue^{\frac{N-2}{N}\gamma}
  \frac{(u+N\hbar)(u-(N-1)\hbar)}
       {(u+\hbar)(u-N\hbar)} \\
& \hspace{0.5 in}\times  \prod_{l=1}^{n-1}\left[
  \frac{(u+N\hbar+N\hbar l)(u-(N-1)\hbar-N\hbar l)}
       {(u+\hbar+N\hbar l)(u-N\hbar-N\hbar l)}\right]e^{-\frac{N-2}{Nl}}, 
\end{align*}
where $\gamma$ is the Euler constant defined by
\[ \gamma=\lim_{n\rightarrow \infty} 
   \left( \sum_{k=1}^n \frac{1}{k}-\log n \right). \]
\end{thm}
We remark that the following formulae hold.
\begin{align*}
(\text{1}) & \exp[ -\sum_{k>0}\frac{1}{2k}f_k^{\ast}(u)_n v^k]=
             \left[ \frac{g^{\ast}(u-v)_n}{g^{\ast}(u)_n} 
             \right]^{\frac{1}{2}}, 
             \hspace{1.1 in} \text{for $\ast=(I),(II)$}\\
(\text{2}) & \lim_{n\rightarrow \infty} g^{\ast}(u)_n=
             \begin{cases}\displaystyle{
             (u-\frac{N-2}{2}\hbar)
            \frac{\Gamma(\frac{1}{2}-\frac{u}{N\hbar})
                  \Gamma(\frac{1}{2}+\frac{u+\hbar}{N\hbar})}
                 {\Gamma(\frac{1}{2}+\frac{u}{N\hbar})
                  \Gamma(\frac{3}{2}-\frac{u+\hbar}{N\hbar})}}
            & \quad \text{for $\ast=(I)$} \\
             \displaystyle{u
            \frac{\Gamma(\frac{u+\hbar}{N\hbar})
                  \Gamma(1-\frac{u}{N\hbar})}
                 {\Gamma(1-\frac{u+\hbar}{N\hbar})
                  \Gamma(1+\frac{u}{N\hbar})}}
            & \quad \text{for $\ast=(II)$}
            \end{cases} \quad .
\end{align*}
The second formula can be proved by using the famous Weierstrass formula
for the Gamma function
\[ \frac{1}{\Gamma(u)}=ue^{\gamma u}
   \prod_{n=1}^{\infty}(1+\frac{u}{n})e^{-\frac{u}{n}}. \]
\begin{rem}
For each $n\in \nz_{>0}$, the fields $\Phi^{(i,i+1)}_{N-1}(u)_n$ and
$\Psi^{(i,i+1)}_{0}(u)_n$ make sense as formal series in $\hbar$. But after
taking the limit $n\rightarrow \infty$, they can not expand with respect to 
$\hbar$. They have to be regarded as, for example, meromorphic functions.
Such feature has never appeared in the quantum affine case \cite{koyama}.
\end{rem}

Here we give a sketch of a proof of Theorem \ref{VOp_sl_N} for type $I$ 
vertex operator and $N=2$ case for simplicity. 
We also give some comments how to prove for general case.

We define the normal ordering $: \cdot :$ of the fields by regarding
$a_{j,k}(k<0), e^{\alpha_j}(1\leq j\leq N-1)$ as creation operators  and
$a_{j,k}(k>0), \partial_{\alpha_j}(1\leq j\leq N-1)$ as annihilation 
operators. After some calculation, we obtain the following operator
product expansion (OPE):
\begin{align*}
\Phi^{(i,i+1)}_{1}(u)_n H_{1}^{-}(v)&=
\left[ \left(\frac{u-v+\frac{3}{4}\hbar}{u-v-\frac{1}{4}\hbar}\right)^2
       \frac{(u-v+(2n-\frac{1}{4})\hbar)(u-v-(2n+\frac{1}{4})\hbar)}
            {(u-v+(2n+\frac{3}{4})\hbar)(u-v-(2n+\frac{5}{4})\hbar)}
       \right]^{\frac{1}{2}} \\
&\times : \Phi^{(i,i+1)}_{1}(u)_n H_{1}^{-}(v) : , \\
\Phi^{(i,i+1)}_{1}(u)_n E_{1}(v)&= (-1)^{\frac{1}{2}}
\left[ (u-v+\frac{1}{2}\hbar)^2
       \frac{u-v-(2n+\frac{1}{2})\hbar}
            {u-v+(2n+\frac{1}{2})\hbar} \right]^{\frac{1}{2}} \\
&\times : \Phi^{(i,i+1)}_{1}(u)_n E_{1}(v) : .
\end{align*}
Taking the limit $n\rightarrow \infty$, fixing the branch, we get
\begin{align*}
& \Phi^{(i,i+1)}_{1}(u)H_{1}^{-}(v)=
   \frac{u-v+\frac{3}{4}\hbar}{u-v-\frac{1}{4}\hbar}
   : \Phi^{(i,i+1)}_{1}(u)H_{1}^{-}(v) :, \\
& \Phi^{(i,i+1)}_{1}(u)E_{1}(v)=
   -(u-v+\frac{1}{2}\hbar): \Phi^{(i,i+1)}_{1}(u)E_{1}(v) :. 
\end{align*}
These are precisely the expected OPE from the intertwining property. 
The other OPE can be obtained easily and here we omit them. Normalization
condition can also be checked similarly.

Next to prove the general case, first simplify OPE, as above, to see the
phase factor and then calculate the limit using the infinite product form
of the Gamma function. In this way, we can prove that our formulae give
the desired OPE and the normalization.

\section{Discussion}\hspace{1 in} \\ 
In this article, we have constructed the Yangian double $\cD Y_{\hbar}(\gtg)$
with a central extension for $\gtg=\gtgl_N. \gtsl_N$. We also presented
Drinfel'd generators which are defined in \cite{dr_2}. Using these 
generators, we studied both finite and infinite dimensional representations.
We presented a conjecture for irreducible finite dimensional representations 
and gave some examples to check the validity of it. The bosonization of 
level $1$ modules and vertex operators were also given.

It seems that the Yangian double $\cD Y_{\hbar}(\gtg)$ for other type of
simple finite dimensional Lie algebra $\gtg$ can be defined by Corollary
\ref{sl_N} without $K^{\pm}(u)$ where $A=(a_{ij})$ is now the corresponding
Cartan matrix. Suppose this is true for a moment. Then the rest of Section
\ref{main} also hold without any change. Especially when $\gtg$ is a simply
laced algebra, we can generalize Theorem \ref{level 1 sl_N} by simple
modification whose quantum affine version is treated in \cite{jing}.
There are several other problems which we have already mentioned in our
previous paper \cite{kohno}. The relation between the quantum affine algebra
$U_q(\hgtg)$ and the Yangian double $\cD Y_{\hbar}(\gtg)$ is quite mysterious.

For physical applications, it is important to investigate the 
infinite dimensional representation theory of $\cD Y_{\hbar}(\gtg)$.
In this paper, we give bosonization of level 1 module $\cF_i$ 
and vertex operators among them. As we have seen in Theorem \ref{VOp_sl_N},
the Fourier coefficients of the vertex operators loose sense 
unlike to the quantum affine case \cite{jimbo_miwa},\cite{koyama}.
This means that we have to consider not the Fourier components but 
the currents themselves. Namely we have to consider a new class of 
the algebra and their representation theory to investigate further.
 It is also interesting to see the connection between the formulae in 
\cite{lukyanov} and ours.

\begin{ack}
The authors would like to thank J. Ding,
L. D. Faddeev, M. Jimbo, M. Kashiwara, M. Kohno, S. Lukyanov, A. Molev, 
T. Miwa, N. Yu. Reshetikhin, E. K. Sklyanin, F. Smirnov and V. Tarasov  
for their interest and valuable suggestions. 
\end{ack}

\appendix
\section{Review of Quantum Groups}\label{quantum_affine}
In this section, we recall some facts about universal $\cR$ and $L$-operator.

\subsection{Universal $\cR$} 
Let $\cR$ be the universal $R$-matrix \cite{dr_icm} for $U_q(\hgtsl_N)$.
For the definition and the properties of universal $R$-matrix, see
\cite{dr_icm,nankai}.

We slightly modify $\cR$ to define $L$-operators. Define
\begin{align*}
& \cR^{'+}=q^{-\frac{1}{2}(c\otimes d+ d\otimes c)}\sigma(\cR^{-1})
           q^{-\frac{1}{2}(c\otimes d+ d\otimes c)}, \\
& \cR^{'-}=q^{\frac{1}{2}(c\otimes d+ d\otimes c)}\cR
           q^{\frac{1}{2}(c\otimes d+ d\otimes c)}, \\
& \cR^{'\pm}(z)=(z^{d}\otimes \id)\cR^{'\pm}(z^{-d}\otimes \id).
\end{align*}
Here $\sigma$ stands for the flip of tensor components
$\sigma(a\otimes b)=b\otimes a$. We remark that $\cR^{'\pm}(z)$ are
formal power series in $z^{\mp 1}$. The properties of universal $R$-matrix
can be readily translated in terms of $\cR^{'\pm}$.
For $x \in U_q(\hgtsl_N)$, we write $\Delta(x)=x_{(1)}\otimes x_{(2)}$.
Then
\[
\cR^{'\pm}(z)\left( \Ad (z^dq^{\pm \frac{1}{2}c_2 d})x_{(1)}\otimes
                    \Ad (   q^{\pm \frac{1}{2}c_1 d})x_{(2)} \right)=
             \left( \Ad (z^dq^{\mp \frac{1}{2}c_2 d})x_{(2)}\otimes
                    \Ad (   q^{\mp \frac{1}{2}c_1 d})x_{(1)} \right)
\cR^{'\pm}(z). 
\]
Here $c_1=c\otimes 1$ and $c_2=1\otimes c$ as in Section \ref{main}.

The Yang-Baxter equation takes the form
\begin{align*}
\cR^{'\pm}_{12}(z/w)\cR^{'\pm}_{13}(zq^{\pm c_2})\cR^{'\pm}_{23}(w)&=
\cR^{'\pm}_{23}(w)\cR^{'\pm}_{13}(zq^{\mp c_2})\cR^{'\pm}_{12}(z/w), \\
\cR^{'+}_{12}(z/wq^{-c_3})\cR^{'+}_{13}(z)\cR^{'-}_{23}(w)&=
\cR^{'-}_{23}(w)\cR^{'+}_{13}(z)\cR^{'+}_{12}(z/wq^{c_3}). 
\end{align*}
For completeness we give the transformation properties of $\cR^{'\pm}$
under the coproduct $\Delta$, the counit $\vep$ and the antipode $S$.
\begin{align*}
& \left( \Delta \otimes \id \right) \cR^{'\pm}(z)=
  \cR^{'\pm}_{13}(zq^{\pm \frac{1}{2}c_2})
  \cR^{'\pm}_{23}(zq^{\mp \frac{1}{2}c_1}), \\
& \left( \id \otimes \Delta \right) \cR^{'\pm}(z)=
  \cR^{'\pm}_{13}(zq^{\mp \frac{1}{2}c_2})
  \cR^{'\pm}_{12}(zq^{\pm \frac{1}{2}c_3}), \\
& (\vep \otimes \id)\cR^{'\pm}(z)=(\id \otimes \vep)\cR^{'\pm}(z)=1, \\
& (S \otimes \id)\cR^{'\pm}(z)=(\id \otimes S)\cR^{'\pm}(z)
 =\cR^{'\pm}(z)^{-1}.
\end{align*}

\subsection{$L$-operators}
Let now $\pi_V:U_q(\hgtsl_N)' \rightarrow \End(V)$ be a finite dimensional 
representation, where $U_q(\hgtsl_N)'$ signifies the subalgebra of
$U_q(\hgtsl_N)$ with $q^d$ being dropped. The evaluation representation 
$\pi_{V_z}$ associated with $V$ is defined by
\[ \pi_{V_z}(x)=\pi_V\left( z^d x z^{-d}\right) \qquad 
   \forall x \in U_q(\hgtsl_N)'.  \]
Introduce the $L$-operators
\begin{align*}
& L^{\pm}(z)=L_{V}^{\pm}(z)=
  \left( \pi_{V_z}\otimes \id \right) \cR^{'\pm}.
\end{align*}
Taking the image of the Yang-Baxter equation for $\cR^{'\pm}$ in
$\End(V_{z})\otimes \End(V_{w})\otimes \id$, we find the following
$RLL$ relations:
\begin{align*}
& R^{\pm}_{12}(z/w)\overset{1}{L^{\pm}}(z)\overset{2}{L^{\pm}}(w)=
  \overset{2}{L^{\pm}}(w)\overset{1}{L^{\pm}}(z)R^{\pm}_{12}(z/w), \\
& R^{+}_{12}(q^{-c}z/w)\overset{1}{L^{+}}(z)\overset{2}{L^{-}}(w)=
  \overset{2}{L^{-}}(w)\overset{1}{L^{+}}(z)R^{+}_{12}(q^{c}z/w),
\end{align*}
where we set
\[ R^{\pm}(z/w)=\left(\pi_{V_{z}}\otimes \pi_{V_{w}} \right) \cR^{'\pm}. \]
Introducing the matrix units $E_{ij}$ let us define the entries 
$L^{\pm}_{ij}(z)$ by
\[ L^{\pm}(z)=\sum E_{ij}\otimes L^{\pm}_{ij}(z). \]
In these terms, the Hopf algebra structure reads as follows.
\begin{align*}
& \Delta \left( L^{\pm}_{ij}(z)\right) =
  \sum_{k} L^{\pm}_{kj}(q^{\pm \frac{1}{2}c_2}z)\otimes
           L^{\pm}_{ik}(q^{\mp \frac{1}{2}c_1}z), \\
& \vep   \left( L^{\pm}_{ij}(z)\right) =
  \delta_{ij}, \\
& S \left({}^t L^{\pm}(z) \right)=\left( {}^t L^{\pm}(z) \right)^{-1}, \\
& S^{-1} \left( L^{\pm}(z)\right)=\left(      L^{\pm}(z) \right)^{-1}.
\end{align*}
In the last two lines we set
\begin{align*}
S \left({}^t L^{\pm}(z) \right) &=
\sum E_{ji} \otimes S \left( L^{\pm}_{ij}(z) \right), \\
S^{-1} \left( L^{\pm}(z)\right) &=
\sum E_{ij} \otimes S^{-1} \left( L^{\pm}_{ij}(z) \right).
\end{align*}
Let $U^{\pm}$ be Hopf subalgebras of $U_q(\hgtsl_N)$ generated by 
$q^{\pm \frac{1}{2}c}$ and the Fourier components of $L^{\pm}(z)$.
The subalgebra $U^{-}$ is the dual Hopf algebra of $U^{+}$ with opposite
comultiplication and the Hopf pairing between $U^{\pm}$ has the explicit
description as follows:
\[
\langle L^{+}(z), L^{-}(w) \rangle =
\sum \langle L^{+}_{ij}(z), L^{-}_{kl}(w) \rangle E_{ij}\otimes E_{kl} =
R^{+}(z/w). 
\]

We remark that all of these formula motivate our choice of 
$T^{\pm}(u)$-matrix.

\section{Several fomulas for $T$-matrix}\label{fusion_T}
In this section, we collect some formulas which seem well-known to 
the specialists \cite{tarasov}. 
 Here we denote $T(u)$ for $T^{\pm}(u)$ for simplicity.

\subsection{Quantum determinant of $T$-matrix}\label{quantum minor}
                                              
In this subsection, we give a brief review on quantum determinant 
for convenience. See \cite{molev},\cite{fusion_text} for further 
information.

\noindent{\bf 1. Quantum minor}

Let $V$ be a rank $N$ $\cA$-free module and $\cP \in \End(V\otimes V)$
be a permutation operator $\cP v\otimes w =w\otimes v~(v,w\in V)$.
Let us fix the normalization of the Yang's $R$-matrix as
\[ R(u)=I+\frac{\hbar}{u}\cP \in \End(V\otimes V). \]
Recall that $T(u)$ enjoy the following commutation relations:
\[ R(u-v)\overset{1}{T}(u)\overset{2}{T}(v)=
   \overset{2}{T}(v)\overset{1}{T}(u)R(u-v). \]
Suppose the comultiplication of $T(u)$ is given by
\[ \Delta(T(u))=T(u)\overset{\cdot}{\otimes}T(u)\quad 
   \text{or equivalently}\quad
   \Delta(t_{ij}(u))=\sum_{k=1}^N t_{ik}(u)\otimes t_{kj}(u). \]
For simplicity, set $R_{i,j}=R_{i,j}(u_i-u_j)$ and
\[ R(u_1,u_2,\cdots ,u_p)=(R_{p-1,p})(R_{p-2,p}R_{p-2,p-1})\cdots
                          (R_{1,p}R_{1,p-1}\cdots R_{1,2}), \]
where the meaning of the lower indices are 
the same as in Section \ref{main}.
\begin{lemma} 
\[ R(u_1,u_2,\cdots ,u_p)
   \overset{1}{T}(u_1)\overset{2}{T}(u_2)\cdots \overset{p}{T}(u_p)
  =\overset{p}{T}(u_p)\cdots \overset{2}{T}(u_2)\overset{1}{T}(u_1)
   R(u_1,u_2,\cdots ,u_p). \]
\end{lemma}

Let $\cA[\gtS_p]$ be the group algebra of the $p$-th symmetric group 
over $\cA$ which naturally acts on $V^{\otimes p}$ and set
\[ a_p=\sum_{\sigma \in \gtS_p}(\sgn \sigma)\sigma \in \cA[\gtS_p],\quad
   A_p=\frac{1}{p!}a_p. \] 

\begin{lemma}[\cite{molev}] For $u_i-u_{i+1}=-\hbar \quad 1\leq \forall i <p$,
\[ R(u_1,u_2,\cdots u_p)=a_p. \]
\end{lemma}
One can prove this lemma by induction on $p$.
Combining these two lemmas, we obtain the following.
\begin{lemma}\label{fused T}
\begin{align*}
 & A_p \overset{1}{T}(u-\frac{p-1}{2}\hbar)
       \overset{2}{T}(u-\frac{p-3}{2}\hbar)\cdots 
       \overset{p}{T}(u+\frac{p-1}{2}\hbar) \\
  = &  \overset{p}{T}(u+\frac{p-1}{2}\hbar)\cdots
       \overset{2}{T}(u-\frac{p-3}{2}\hbar)
       \overset{1}{T}(u-\frac{p-1}{2}\hbar) A_p. 
\end{align*}
\end{lemma}
Set $p=N$ in the above lemma.
Since $N$-th. exterior power $\bigwedge^N V$ is of rank $1$ and $A_N$
stabilizes $\bigwedge^N V$, the left hand side of the above equation
is $(\text{scalar})\times A_N$.
\begin{defn}[Quantum Determinant]
\[ \qdet. T(u)A_N=
       \overset{N}{T}(u+\frac{N-1}{2}\hbar)\cdots
       \overset{2}{T}(u-\frac{N-3}{2}\hbar)
       \overset{1}{T}(u-\frac{N-1}{2}\hbar) A_N. \]
\end{defn}
Explicitly, we have
\begin{prop}\label{defn_q-det}
\begin{align*}
\qdet T(u) & =\sum_{\sigma \in \gtS_N}(\sgn \sigma)
              t_{\sigma(1),1}(u-\frac{N-1}{2}\hbar)
              t_{\sigma(2),2}(u-\frac{N-3}{2}\hbar)\cdots
              t_{\sigma(N),N}(u+\frac{N-1}{2}\hbar) \\
           & =\sum_{\sigma \in \gtS_N}(\sgn \sigma)
              t_{1,\sigma(1)}(u+\frac{N-1}{2}\hbar)
              t_{2,\sigma(2)}(u+\frac{N-3}{2}\hbar)\cdots
              t_{N,\sigma(N)}(u-\frac{N-1}{2}\hbar). 
\end{align*}
\end{prop}

Next we explain some facts about quantum minors of $T$-matrix.
For two index subsets $I,J \subset \{ 1,2 \cdots N \}$ with 
$\#I=\#J=p~,~1\leq p\leq N$ (the cardinality), set 
\[ T_{IJ}(u)=(t_{ij}(u))_{i\in I,j\in J}. \]
By definition of $T(u)$, we obtain the following commutation relations
\[ R_p(u-v)\overset{1}{T}_{IJ}(u)\overset{2}{T}_{IJ}(v)
  =\overset{2}{T}_{IJ}(v)\overset{1}{T}_{IJ}(u)R_p(u-v),\quad
   R_p(u)=I+\frac{\hbar}{u}\cP \in \End(V_p\otimes V_p), \]
where $V_p$ is a rank $p$ $\cA$-free module. Thus by the similar argument
as above, we get the explicit expression of quantum minor $\qdet T_{IJ}(u)$
as follows. Set
\[ I=\{ i_1,i_2,\cdots ,i_p \},J=\{j_1,j_2,\cdots ,j_p \}. \]
\begin{lemma}\label{defn_q-minor}
\begin{align*}
\qdet T_{IJ}(u) & =\sum_{\sigma \in \gtS_p}(\sgn \sigma)
              t_{i_{\sigma(1)},j_1}(u-\frac{p-1}{2}\hbar)
              t_{i_{\sigma(2)},j_2}(u-\frac{p-3}{2}\hbar)\cdots
              t_{i_{\sigma(N)},j_N}(u+\frac{p-1}{2}\hbar) \\
                & =\sum_{\sigma \in \gtS_p}(\sgn \sigma)
              t_{i_1,j_{\sigma(1)}}(u+\frac{p-1}{2}\hbar)
              t_{i_2,j_{\sigma(2)}}(u+\frac{p-3}{2}\hbar)\cdots
              t_{i_N,j_{\sigma(N)}}(u-\frac{p-1}{2}\hbar).
\end{align*}
\end{lemma}
The following Corollary is the immediate consequence of the above expression.
\begin{cor}\label{prop_q-minor} For each $\sigma \in \gtS_p$,
\[ \qdet T_{I^{\sigma}J}(u)=\qdet T_{IJ^{\sigma}}(u)
  =(\sgn \sigma)\qdet T_{IJ}(u), \]
where we set
$I^{\sigma}=\{ i_{\sigma^{-1}(1)},i_{\sigma^{-1}(2)},\cdots ,
               i_{\sigma^{-1}(p)} \}$ 
and similarly for $J^{\sigma}$.
\end{cor}
Using this Corollary, one can calculate the coproduct of quantum minors as
follows.
\begin{cor}\label{coprod_q-minor}
\[ \Delta(\qdet T_{IJ}(u))=\sum_{K}\qdet T_{IK}(u)\otimes \qdet T_{KJ}(u), \]
where the summation runs over all of the ordered subset 
$K=\{k_1,k_2,\cdots, k_p\} \subset S$ satisfying 
$1\leq k_1<k_2<\cdots <k_p\leq N$.
\end{cor}

\noindent{\bf 2. Laplace expansion of $T$-matrix}

Let $\{ e_i \}_{1\leq i \leq N}$ be an $\cA$-free basis of $V$ and 
$S=\{1,2, \cdots ,N \}$ be the index set. For each 
ordered index subset $I=\{ i_1,i_2,\cdots ,i_p \} \subset S$, we define
$e_I$ an element of $\bigwedge^p V$ as 
\[ e_I=\sum_{\sigma \in \gtS_p}(\sgn \sigma)
          e_{i_{\sigma(1)}}\otimes e_{i_{\sigma(2)}}\otimes \cdots \otimes
          e_{i_{\sigma(p)}}. \]
Note that the set $\{ e_I \}_{1\leq i_1 < i_2 < \cdots < i_p \leq N}$ provides
a basis of $\bigwedge^p V$. Let $E_{IJ}$ be an element of 
$\End(\bigwedge^p V)$ satisfying $E_{IJ}.e_{K}=\delta_{JK}e_I$. Set
\[ T_p(u)= \overset{p}{T}(u+\frac{p-1}{2}\hbar)\cdots
           \overset{2}{T}(u-\frac{p-3}{2}\hbar)
           \overset{1}{T}(u-\frac{p-1}{2}\hbar) A_p. \]
By Lemma \ref{fused T}, we see that $T_p(u)$ is the element of 
$\End(\bigwedge^p V)$. More precisely, we have
\begin{lemma}\label{key_T}
\[ T_p(u)e_J=\sum_{I} (\qdet T_{IJ}(u))e_{I} \quad \text{or equivalently}
   \quad T_p(u)=\sum_{I,J}(\qdet T_{IJ}(u))E_{IJ}. \]
\end{lemma}
One can prove this lemma by using Corollary \ref{prop_q-minor}.

Fix $p,q \in \nz_{>0}$ such that $p+q=N$. Regarding both $\bigwedge^N V$ and
$\bigwedge^p V \otimes \bigwedge^q V$ as subspaces of $V^{\otimes N}$, one can 
easily express $e_S \in \bigwedge^N V$ by linear combinations of 
$e_I\otimes e_J \in \bigwedge^p V \otimes \bigwedge^q V$ as follows.

\begin{lemma}\label{e_S}
\[ e_S=\sum_{I\cup J=S;\#I=p}(-1)^{\frac{1}{2}p(p+1)+|I|}e_I\otimes e_J, \]
where $|I|=\sum_{j=1}^p i_j$ for $I=\{i_1,i_2,\cdots ,i_p \}$.
\end{lemma}

Combining Lemma \ref{key_T} and Lemma \ref{e_S}, we obtain the Laplace 
expansion of $T$-matrix.
\begin{prop}[Quantum Laplace Expansion]\label{Laplace_expansion}
For each $I,J \subset S$ with $\#I=p,\#J=q$, we have
\[ (\qdet T(u))\delta_{I\cup J,S}\delta_{I\cap J,\phi}
=\sum_{K\cup L=S; \# K=p}(-1)^{|I|+|K|}
 \qdet T_{IK}(u+\frac{q}{2}\hbar)\qdet T_{JL}(u-\frac{p}{2}\hbar). \]
\end{prop}
Specializing $p=1$ or $q=1$ we obtain quantum minor expansion of $T$-matrix.
Namely set $S^{(i)}=S\setminus \{i\}$ and
\[ \widetilde{T}(u)=(\tilde{t}_{ij}(u))_{1\leq i,j \leq N},\quad
   \tilde{t}_{ij}(u)=(-1)^{i+j}\qdet T_{S^{(j)},S^{(i)}}(u). \]
Then we get
\begin{cor}\label{inverse of T}
\[ T(u+\frac{N-1}{2}\hbar)\widetilde{T}(u-\frac{1}{2}\hbar)
  ={}^{t}\widetilde{T}(u+\frac{1}{2}\hbar){}^{t}T(u-\frac{N-1}{2}\hbar)
  =(\qdet T(u))I, \]
where ${}^{t}$ signifies the transpose of the matrix.
\end{cor}

\subsection{Gauss decomposition of $T$-matrix}\label{gauss}
In this subsection, we explicitly construct Gauss decomposition of
$T(u)$ in terms of their quantum minors.
Let  
\begin{align*}
T(u) &= \begin{pmatrix}
        1 & & & 0 \\
        f_{2,1}(u) & \ddots & & \\
            &\ddots &\ddots & \\
        f_{N,1}(u)& &f_{N,N-1}(u) & 1
          \end{pmatrix}
      \begin{pmatrix}
      k_{1}(u) & & & 0 \\
              & \ddots & &   \\
              & & \ddots & \\
      0 & & & k_{N}(u)
              \end{pmatrix} \\
          & \times \begin{pmatrix}
            1 & e_{1,2}(u)& & e_{1,N}(u) \\
              & \ddots &\ddots & \\
              & & \ddots & e_{N-1,N}(u) \\
            0 & & & 1
            \end{pmatrix}.
\end{align*}
be the Gauss decomposition of $T(u)=(t_{ij}(u))$. 
\begin{lemma}\label{gauss_comp}
\[ t_{i,j}(u)=\begin{cases}
   \sum_{l<i}f_{i,l}(u)k_{l}(u)e_{l,j}(u)+k_{i}(u)e_{i,j}(u)
   &  i<j , \\
   \sum_{l<i}f_{i,l}(u)k_{l}(u)e_{l,i}(u)+k_{i}(u)
   &  i=j , \\
   \sum_{l<j}f_{i,l}(u)k_{l}(u)e_{l,j}(u)+f_{i,j}(u)k_{i}(u)
   &  i>j .   \end{cases} \]
\end{lemma}

For $1\leq p,q \leq N$, let us define $T_{p,q}(u)$ submatrices of
$T(u)$ as follows.
\begin{defn}
\begin{align*}
(1)& p=q, \qquad T_{p,p}(u)=
          \left( t_{i,j}(u) \right)_{1\leq i,j \leq p} \\
(2)& p<q, \\
   & T_{p,q}(u)=\begin{pmatrix}
     t_{1,1}(u)  & \dots & t_{1,p-1}(u)  & t_{1,q}(u) \\
     \vdots      &       & \vdots        & \vdots     \\
     t_{p-1,1}(u)& \dots & t_{p-1,p-1}(u)& t_{p-1,q}(u) \\
     t_{p,1}(u)  & \dots & t_{p,p-1}(u)  & t_{p,q}(u)
                \end{pmatrix} , \\
(3)& p>q, \\
   & T_{p,q}(u)=\begin{pmatrix}
     t_{1,1}(u)  & \dots & t_{1,q-1}(u)  & t_{1,q}(u) \\
     \vdots      &       & \vdots        & \vdots     \\
     t_{q-1,1}(u)& \dots & t_{q-1,q-1}(u)& t_{q-1,q}(u) \\
     t_{p,1}(u)  & \dots & t_{p,q-1}(u)  & t_{p,q}(u)
                \end{pmatrix} . 
\end{align*}
\end{defn}
Using Lemma \ref{gauss_comp}, we can explicitly describe the 
Gauss decomposition of $T_{p,q}(u)$ as follows.
\begin{lemma}\label{gauss_submatrix}
\begin{align*}
(1) & p=q, \\
 T_{p,p}(u)&=\begin{pmatrix}
      1         &        &             & 0       \\
      f_{2,1}(u)& \ddots &             &         \\
      \vdots    &        & \ddots      &         \\
      f_{p,1}(u)& \hdots & f_{p,p-1}(u)& 1
                 \end{pmatrix}
                 \begin{pmatrix}
      k_{1}(u)  &        &             & 0       \\
                & \ddots &             &         \\
                &        & \ddots      &         \\
          0     &        &             & k_{p}(u)
                 \end{pmatrix}  \\
    & \times \begin{pmatrix}
      1         & e_{1,2}(u)& \hdots   & e_{1,p}(u)\\
                & \ddots    & \ddots   & \vdots    \\
                &           & \ddots   & e_{p-1,p}(u) \\
          0     &           &          &    1      
                 \end{pmatrix} \\
(2) & p<q, \\
T_{p,q}(u)&=\begin{pmatrix}
      1           &        &               & 0       & 0 \\
      f_{2,1}(u)  & \ddots &               &         & \\
      \vdots      & \ddots & \ddots        &         & \\
      f_{p-1,1}(u)& \hdots & f_{p-1,p-2}(u)& 1       & \\
      f_{p,1}(u)  & \hdotsfor{2}       & f_{p,p-1}(u)& 1 
                 \end{pmatrix}
                 \begin{pmatrix}
      k_{1}(u)    &        &               &         & 0 \\
                  & \ddots &               &         &   \\
                  &        & \ddots        &         &   \\
                  &        &             & k_{p-1}(u)&   \\
         0        &        &      &  & k_{p}(u)e_{p,q}(u)
                 \end{pmatrix} \\
    & \times \begin{pmatrix}
      1           &e_{1,2}(u)& \hdots & e_{1,p-1}(u)& e_{1,q}(u) \\
                  & \ddots   & \ddots & \vdots      & \vdots     \\
                  &          & \ddots & e_{p-2,p-1}(u)& \vdots   \\
      0           &          &        &     1       & e_{p-1,q}(u) \\
      0           &          &        &             &    1        
                 \end{pmatrix} \\
(3) & p>q, \\
T_{p,q}(u)&=\begin{pmatrix}
            1           &        &               & 0       & 0 \\
      f_{2,1}(u)  & \ddots &               &         & \\
      \vdots      & \ddots & \ddots        &       & & \\
      f_{q-1,1}(u)& \hdots & f_{q-1,q-2}(u)&  1    & \\
      f_{p,1}(u)  & \hdotsfor{2}       & f_{p,q-1}(u)& 1 
                 \end{pmatrix}
                 \begin{pmatrix}
      k_{1}(u)    &        &               &         & 0 \\
                  & \ddots &               &         &   \\
                  &        & \ddots        &         &   \\
                  &        &             & k_{q-1}(u)&   \\
         0        &        &      &  & f_{p,q}(u)k_{q}(u)
                 \end{pmatrix} \\
    & \times \begin{pmatrix}
      1           &e_{1,2}(u)& \hdots & e_{1,q-1}(u)& e_{1,q}(u) \\
                  & \ddots   & \ddots & \vdots      & \vdots     \\
                  &          & \ddots & e_{q-2,q-1}(u)& \vdots   \\
      0           &          &        &     1       & e_{q-1,q}(u) \\
      0           &          &        &             &    1        
                 \end{pmatrix} 
\end{align*}
\end{lemma}
Let
\[ T_{p,q}(u)=F_{p,q}(u)K_{p,q}(u)E_{p,q}(u) \]
be the Gauss decomposition of $T_{p,q}(u)$ and set $r=\min\{p,q\}$.
Comparing the $(r,r)$ component of the formula
\[ F_{p,q}(u)^{-1}=K_{p,q}(u)E_{p,q}(u)T_{p,q}(u)^{-1} \]
in the both hand side together with Lemma \ref{gauss_submatrix},
we obtain the following.
\begin{lemma}\label{gauss_pre}
\begin{align*}
(1) & k_{p}(u)=\dfrac{1}{\left[ T_{p,p}(u)^{-1} \right]_{p,p}}, \\
(2) & e_{p,q}(u)=\left[ T_{p,p}(u)^{-1} \right]_{p,p} 
               \dfrac{1}{\left[ T_{p,q}(u)^{-1} \right]_{p,p}}, \\
(3) & f_{p,q}(u)=\dfrac{1}{\left[ T_{p,q}(u)^{-1} \right]_{q,q}}
                 \left[ T_{q,q}(u)^{-1} \right]_{q,q}, 
\end{align*}
where $\left[ T_{p,q}(u)^{-1} \right]_{a,b} $ signifies the $(a,b)$
component of the matrix $T_{p,q}(u)^{-1}$.
\end{lemma}

Set \[
\Delta_{p,q}(u):=\qdet. T_{p,q}(u),\qquad 
\Delta_{p}(u):=\qdet. T_{p,p}(u). \]
Since we can express matrix components of $T_{p,q}(u)^{-1}$ by
their quantum minors Lemma \ref{inverse of T}, combining with Lemma
\ref{gauss_pre}, we obtain the following results.
\begin{thm}\label{gauss_decomp}
\begin{align*}
(1) & k_{p}(u)=\Delta_{p}\left(u-\frac{p-1}{2}\hbar \right)
               \Delta_{p-1}\left(u-\frac{p}{2}\hbar \right)^{-1}, \\
(2) & e_{p,q}(u)=\Delta_{p}\left(u-\frac{p-1}{2}\hbar \right)^{-1}
               \Delta_{p,q}\left(u-\frac{p-1}{2}\hbar \right), \\
(3) & f_{p,q}(u)=\Delta_{p,q}\left(u-\frac{q-1}{2}\hbar \right)
               \Delta_{q}\left(u-\frac{q-1}{2}\hbar \right)^{-1}.
\end{align*}
\end{thm}

\end{document}